\documentclass[11pt,sort&compress]{elsarticle}
\usepackage{amssymb,epsfig,multirow,caption,subcaption,calculator}
\usepackage{natbib}
\usepackage{graphicx}
\usepackage{array}
\usepackage{amsmath}
\usepackage{amsfonts}
\usepackage{wrapfig}
\usepackage[font=small,labelfont=bf]{caption}
 \usepackage{tikz}
\usepackage{tabu}
\usepackage{floatrow}
\newcommand{\dt}{\Delta t}

\begin{document}

\begin{frontmatter}
\title{Adaptive pseudo-time methods for the Poisson-Boltzmann equation with Eulerian solvent excluded surface}
\author[label1]{Benjamin Jones}
\address[label1]{Department of Mathematics, University of Alabama, Tuscaloosa, AL 35487,USA}
\author[label2]{Sheik Ahmed Ullah}
\address[label2]{Department of Mathematics, Stillman College, Tuscaloosa, AL 35401,USA}
\author[label1]{Siwen Wang}
\cortext[cor1]{Corresponding author}
\author[label1]{Shan Zhao\corref{cor1}}
\ead{szhao@ua.edu}

\begin{abstract}
This work further improves the pseudo-transient approach for the Poisson Boltzmann equation (PBE) in the electrostatic analysis of solvated biomolecules. The numerical solution of the nonlinear PBE is known to involve many difficulties, such as exponential nonlinear term, strong singularity by the source terms, and complex dielectric interface. Recently, a pseudo-time ghost-fluid method (GFM) has been developed in [S. Ahmed Ullah and S. Zhao, {\it Applied Mathematics and Computation}, {\bf 380}, 125267, (2020)], by analytically handling both nonlinearity and singular sources. The GFM interface treatment not only captures the discontinuity in the regularized potential and its flux across the molecular surface, but also guarantees the stability and efficiency of the time integration. However, the molecular surface definition based on the MSMS package is known to induce instability in some cases, and a nontrivial Lagrangian-to-Eulerian conversion is indispensable for the GFM finite difference discretization. In this paper, an Eulerian Solvent Excluded Surface (ESES) is implemented to replace the MSMS for defining the dielectric interface. The electrostatic analysis shows that the ESES free energy is more accurate than that of the MSMS, while being free of instability issues. Moreover, this work explores, for the first time in the PBE literature, adaptive time integration techniques for the pseudo-transient simulations. A major finding is that the time increment $\dt$ should become smaller as the time increases, in order to maintain the temporal accuracy. This is opposite to the common practice for the steady state convergence, and is believed to be due to the PBE nonlinearity and its time splitting treatment. Effective adaptive schemes have been constructed so that the pseudo-time GFM methods become more efficient than the constant $\dt$ ones. 
\end{abstract}

\end{frontmatter}


\section{Introduction}
The Poisson Boltzmann equation (PBE)  \cite{Honig95} is a widely used implicit solvent continuum model for studying electrostatic interaction between a biomolecule and its solvent environment. Such an electrostatic analysis plays an important role in understanding the structure, function, and dynamics of the biomolecules.  By treating the solute biomolecule and the solvent, respectively, as low-dielectric and high-dielectric medium, the PBE is formulated as a nonlinear elliptic partial differential equation  with singular source terms to account for partial charges contained in the macromolecule. The numerical solution of the three-dimension (3D) PBE for real protein structures is known to be challenging in various aspects, such as strong singularity due to the charge sources \cite{Geng17}, geometrically complicated molecular surface as the solute-solvent boundary \cite{Connolly_surface,Sanner_MSMS}, the exponential nonlinearity for strong ionic strength \cite{Geng13,Zhao14}, and special accommodation needed for the regularity loss of the potential across the dielectric interface \cite{Chen11}. To deal with these difficulties, various software packages have been successfully developed in the literature for electrostatic analysis, such as DelPhi \cite{Li:2012}, AMBER \cite{AMBER}, APBS \cite{APBS}, MIBPB \cite{Chen11}, and SDPBS \cite{Xie:2014}.

Recently, a pseudo-transient approach \cite{sheik_zhao_2020} has been introduced for solving the PBE, which is equipped with effective components to overcome the aforementioned numerical difficulties. In treating singular charge sources, the two-component regularization method developed in \cite{Geng17} is adopted so that the singular component of the potential can be analytically captured by Green's functions. One then needs to solve a regularized PBE with both the solution and its flux being discontinuous across the dielectric interface. To deal with the hyperbolic sine term, the regularized PBE boundary value problem is converted into a pseudo-time steady-state problem. In a time splitting setting, the nonlinear subsystem can be analytically integrated \cite{Geng13,Zhao14} so that the nonlinearity instability is suppressed. The time stepping of the linear subsystem can be carried out by using fully implicit alternating direction implicit (ADI) \cite{Geng13,Zhao14} or locally one dimensional (LOD)  \cite{Leighton_Zhao} methods, both of which efficiently reduce the 3D problems into one-dimensional (1D) ones. Moreover, the ghost-fluid method (GFM) originally developed in  \cite{Fedkiw1999,Liu2000} has been reformulated in the ADI framework \cite{sheik_zhao_2020} for treating complicated molecular surface and discontinuous solution across the dielectric interface. The modified GFM minimizes the geometric information needed for the molecular surface, and maintains the symmetry and diagonal dominance of the 1D finite-difference matrix, so that the stability of ADI and LOD methods is boosted \cite{sheik_zhao_2020}. 

The first goal of this work is to further improve the pseudo-time GFM approach by adopting a more suitable molecular surface. 
In the PBE literature, the solute-solvent boundary is usually taken as the solvent excluded surface (SES) \cite{Lee:1971,Richards:1977}, which is defined with a solvent molecule probe that rolls around the Van der Waals spheres representing all atoms. The boundary of the surface is determined by locations that the probe cannot access. Connolly described how this surface can be analytically calculated through a series of patches \cite{Connolly_surface}. The MSMS software developed in \cite{Sanner_MSMS} provides an efficient generation of the SES based on a reduced surface, and has been adopted in many PBE solvers, as well as in visualization softwares such as Chimera \cite{chimera} and VMD \cite{VMD}. Mathematically, the MSMS output provides a Lagrangian representation of the SES by means of a triangulated surface mesh, whose quality is controlled by a specified density measuring triangle vertices per angstrom \cite{Sanner_MSMS}. In the pseudo-time GFM approach \cite{sheik_zhao_2020}, in order to facilitate the GFM interface treatments, a Lagrangian-to-Eulerian algorithm similar to the one presented in \cite{Yu2007} has been developed to convert the MSMS triangulation into Cartesian grids. 

However, some problems have been reported for the MSMS in the literature.  
For example,  if an atom causes MSMS to fail to compute a surface, MSMS will increase the atomic radius of that atom by $0.1$\AA~and restart the surface calculation \cite{Sanner_MSMS}. After three failures, MSMS quits and fails to generate a surface, which happens in several large systems \cite{TMSMesh}.
In other cases, MSMS generates an incorrect surface \cite{nanoshaper,ESES}, especially when large densities are employed. Moreover, the Lagrangian-to-Eulerian conversion is numerically challenging, because the SES molecular surface is known to have geometrical singularities \cite{Yu2007}. Because the GFM scheme only needs the interface locations, and not the normal directions, the conversion in the pseudo-time approach \cite{sheik_zhao_2020} is slightly simpler than that of the MIB method \cite{Yu2007}. Still, the computation becomes nontrivial when the density is large, because with smaller sizes of triangles, the interaction of one Cartesian grid line with the SES is hard to find. The MSMS induced instability is frequently encountered in the PBE simulations. 

Fortunately, an Eulerian Solvent Excluded Surface (ESES) algorithm \cite{ESES} has been constructed recently to directly calculate the analytical SES patches given by Connolly \cite{Connolly_surface} on Cartesian grids. The ESES computes distances between Cartesian grid points and surface intersections as well as the normal direction of the surface at the intersection, which can be directly adopted in the pseudo-time GFM approach without further conversion. Moreover, when applying the ESES in the existing MIB-PBE solver \cite{Chen11}, it has been observed  \cite{ESES,ESES_parallel} that the ESES surface is not only free of instability issues, but also produces more accurate free energies. In particular, while being stable, the MSMS solvation energy usually approaches to that of the ESES as the MSMS density increases. This motivates us to replace the MSMS by the ESES in the pseudo-time GFM approach and investigate the stability and performance of the improved pseudo-time solver in this study. 

The second goal of this work is to improve the computational efficiency of the pseudo-time GFM approach. In the pseudo-transient approach, the solution of the PBE is recovered from the steady-state solution of the pseudo-time dependent PBE \cite{Geng13,Zhao14}, which means that a long time integration is required. In order to allow the use of a large time increment $\dt$ for efficiency, implicit methods are normally employed in pseudo-time approaches \cite{Geng13,Leighton_Zhao,DENG2018}. For implicit schemes, a large linear system obtained from discretizing the 3D PBE has to be solved in each time step. Consider the spatial degree of freedom to be $N$ for the 3D system. Both alternating direction implicit (ADI)  \cite{Geng13,Zhao14} and locally one dimensional (LOD)  \cite{Leighton_Zhao} methods have been developed to achieve the fastest speed for solving the PBE linear systems. By converting 3D systems into a set of independent 1D systems, the complexity of ADI and LOD is just $O(N \log N)$ for each time step. 

In the present study, we will develop an adaptive pseudo-time approach for the PBE. Our aim is to select $\dt$ adaptively so that total time steps could be minimized, without sacrificing the temporal accuracy. In the literature, the use of adaptive time step is a well known technique for an efficient convergence to the steady state \cite{Bank1980,Kelley1998,Pollock2015}. Because the solution approaches to certain limit as the time $t$ becomes larger, the temporal variation of the solution becomes smaller. Thus, one tends to use a large $\dt$ as $t$ increases to save total time steps. However, the present study shows contradictory results. For the pseudo-time GFM approach, the final accuracy in estimating the electrostatic free energy depends on $\dt$ critically. In order to maintain such accuracy, a small enough $\dt$ is required before the steady state. Therefore, our strategy for minimizing the time steps is using a large $\dt$ initially, and reduce it as $t$ increases. Moreover, the proportional–integral–derivative (PID) method \cite{pid} originally developed for fluid dynamics is employed to select $\dt$ adaptively. 

The rest of this paper is organized as follows. In the Section 2, we will briefly review the PBE model and the pseudo-time GFM approach. The replacement of the MSMS by the ESES will be discussed in the Section 3, together with numerical validation of the new pseudo-time GFM solver based on the ESES. The adaptive time selection of this new solver will be considered in the Section 4. Several strategies will be explored so that a tradeoff could be achieved in minimizing the CPU time and preserving the accuracy. A large scale tests of proteins will be carried out to benchmark the new PBE solvers.  Finally, this paper ends with a conclusion.

\section{Poisson-Boltzmann equation and pseudo-time methods}
In this section, we briefly review the physical model and the pseudo-time ghost-fluid method (GFM) developed in \cite{sheik_zhao_2020}. The proposed numerical improvements will be presented in next two sections. 

\subsection{Poisson-Boltzmann Equation}
\begin{figure}
    \centering
    \includegraphics[width=.6\textwidth]{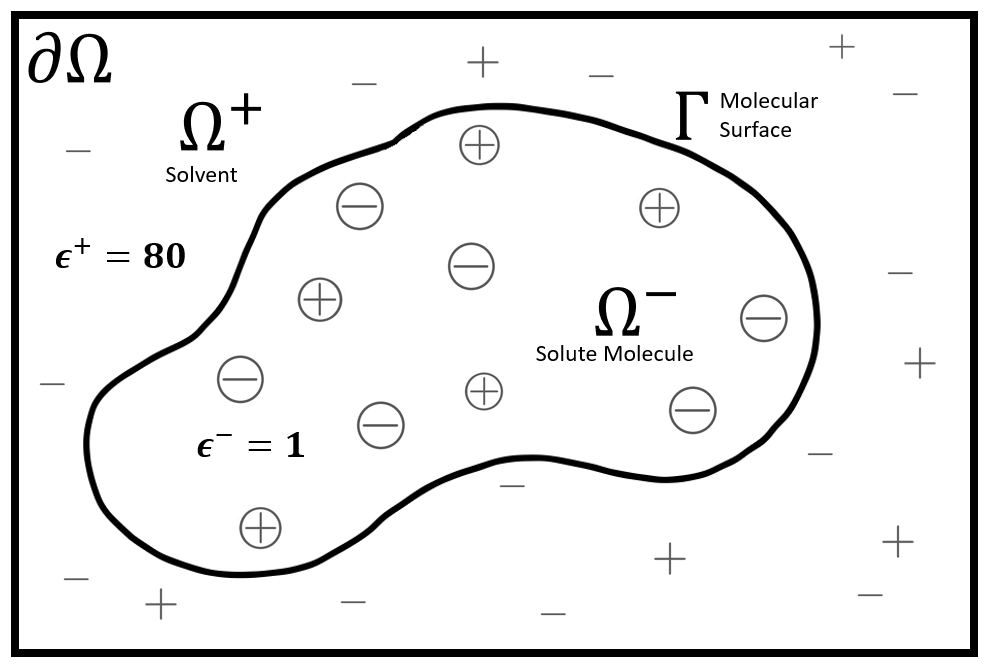}
    \caption{Illustration of the PBE domain.}
    \label{fig:domain}
\end{figure}

The Poisson-Boltzmann Equation (PBE) is the governing equation of electrostatics for a solute macromolecule immersed in an aqueous solvent environment \cite{Honig95}. After loading a protein structure from protein databank \cite{RCSB_PDB}, a large enough cubic domain $\Omega \in \mathbb{R}^3$ is first identified as the computational domain. This domain consists of two regions, $\Omega^-$ and $\Omega^+$, respectively, for solute and solvent, with the solute-solvent boundary defined by the molecular surface $\Gamma$ \cite{Lee:1971,Richards:1977}. A two-dimensional representation of this domain can be seen in Fig. \ref{fig:domain}.   The electrostatic interaction of this solute-solvent system for $\textbf{r} \in \Omega$ is governed by the nonlinear PBE as 
\begin{equation}\label{pbe}
			-\nabla \cdot(\epsilon(\textbf{r})\nabla \phi(\textbf{r}))+\kappa^2(\textbf{r}) \sinh (\phi(\textbf{r}))=\rho(\textbf{r}), 
\end{equation}
where $\phi(\textbf{r})$ is the electrostatic potential and the singular source $\rho(\textbf{r})$ term is defined as
\begin{equation}\label{rho}
	\rho(\textbf{r})= 4\pi \frac{e_c^2}{k_B T}\sum_{i=1}^{N_c} q_i \delta(\textbf{r}-\textbf{r}_i). 
\end{equation}

Here $N_c$ is the total number of atoms in the solute molecule, $T$ is the temperature, $k_B$ is the Boltzmann constant, $e_c$ is the fundamental charge and $q_i$, in the same unit as $e_c$ is the partial charge on the \textit{i}th atom of the solute molecule located at position $\textbf{r}_i$. 
On the outer boundary of $\Omega$, i.e., $\partial \Omega$, a Dirichlet boundary condition is usually assumed in biomolecular simulations 
\begin{equation*}
	\phi (\textbf{r})=\phi_b (\textbf{r}) := \frac{e_c^2}{k_B T} \sum_{i=1}^{N_c} \frac{q_i e^{-|\textbf{r}-\textbf{r}_i | \sqrt{\frac{\kappa^2}{\epsilon^+}} }}{\epsilon^{+}|\textbf{r}-\textbf{r}_i|}, \quad \text{ on } \partial \Omega.
\end{equation*}

The PBE is known as an elliptic interface problem in the literature, because across the interface $\Gamma$, two parameters presented in (\ref{pbe}) are piecewisely defined. The relative permittivity $\epsilon$ takes a low dielectric value $\epsilon^-$ in $\Omega^-$ and a high dielectric value $\epsilon^+$ in $\Omega^+$. The Debye-Huckel parameter $\kappa^2 = 8.486902807$\AA$^{-2} I_s$  \cite{sheik_zhao_2020} for $\textbf{r} \in \Omega^+$ and $\kappa=0$ for $\textbf{r} \in \Omega^-$. Consequently, the potential losses its regularity across $\Gamma$. Numerically, in order to guarantee the order of convergence, two physical conditions for the potential and flux density have to be satisfied in the discretization. These jump conditions are 
\begin{equation}
\left[\phi \right]_\Gamma = 0 \textnormal{ and } \left[\epsilon\phi_n\right]_\Gamma = 0. \label{ju_cond}
\end{equation}
where $\textbf{n} =(n_x,n_y,n_z)$ is the outer normal direction on the interface $\Gamma$ and $\phi_n= \frac{\partial \phi}{\partial n} $ is the directional derivative in $\textbf{n}$. The notation $[f]_\Gamma = f^+-f^-$ represents the difference of the functional value across the interface $\Gamma$. 

\subsection{Pseudo-transient approach}
In the pseudo-transient approach \cite{sheik_zhao_2020}, to overcome the numerical challenge associated with the singular source $\rho$ in (\ref{pbe}), a two-component regularization formulation proposed in \cite{Geng17} is first applied to transform the PBE into a source free regularized PBE. Then in order to bypass the nonlinear instability due to the exponential nonlinear term $\sinh (\phi(\textbf{r}))$, the PBE is converted into a time-dependent steady state problem \cite{Geng13,Zhao14}. The governing system can be given as \cite{sheik_zhao_2020}
\begin{eqnarray}
	 	\frac{\partial u(\textbf{r},t)}{\partial t} &=& \nabla \cdot(\epsilon({\bf r}) \nabla u(\textbf{r},t))-\kappa^2({\bf r}) \sinh(u(\textbf{r},t)),\text{ in } \Omega^- \cup \Omega^+ \label{eq:TRPBE}, \\
	 	  \left[u\right]&=& G, \quad\text{ on } \Gamma,\label{eq:TRPBE_jump1}\\
	 	  \left[\epsilon\frac{\partial u}{\partial n}\right]&=& \epsilon^-  \frac{\partial G}{\partial n}, \quad \text{ on } \Gamma,\label{eq:TRPBE_jump2}\\ 
	 	 u&=&\phi_b, \quad \text{ on } \partial \Omega, \label{eq:bc}
\end{eqnarray}
where $u(\textbf{r},t)$ is the pseudo-time dependent potential. Its steady state solution gives rise to the regularized potential, from which the original potential $\phi$ can be directly calculated \cite{Geng17}. We note that because of regularization, the jump conditions (\ref{eq:TRPBE_jump1}) and (\ref{eq:TRPBE_jump2}) becomes nonhomogeneous, i.e., both the potential and its flux are now discontinuous. Here $G({\bf r})$ is the Green's function due to the singular charges \cite{Geng17} 
\begin{equation}\label{Green} 
	G({\bf r}) = \frac{e_c^2}{k_B T } \sum_{i=1}^{N_c} \frac{q_i }{\epsilon^{-}|{\bf r}-{\bf r}_i|}. 
\end{equation}
Both $G({\bf r})$ and its gradient are analytically defined on $\Gamma$. 

Two efficient operator splitting schemes have been proposed in the PBE literature for analytically treating the nonlinear term $\sinh (u)$, i.e., the alternating direction implicit (ADI) \cite{Geng13,Zhao14} and locally-one-dimensional (LOD) \cite{Leighton_Zhao}. 
Let us consider a uniform mesh with a grid spacing $h$ in all $x,y$ and $z$ directions having $N_x,N_y$ and $N_z$ as the number of the grid points in each direction. We assume the vector $U^n = \left\{u^n_{ijk}\right\}$ for $i=1,...., N_x,j=1,...N_y,$ and $k=1,...N_z$ having all the nodal values of $u$ at the time level $t_n$ as its elements. To update $U^n$ at time level $t_n$ to $U^{n+1}$ at time level $t_{n+1}=t_n+ \Delta t $, a time splitting with two stages is employed in the ADI scheme \cite{Geng13,Zhao14} 
\begin{eqnarray}
  \frac{\partial w}{\partial t}=& -\kappa^2 \sinh(w) &\text{ with } W^n=U^n\text{ and } t \in \left[t_n,t_{n+1}\right],\label{non_linear_ADI}\\
 \frac{\partial v}{\partial t}=&  \nabla \cdot(\epsilon\nabla v) &\text{ with } V^n=W^{n+1}\text{ and } t \in \left[t_n,t_{n+1}\right]	 .\label{diffusion_ADI}
\end{eqnarray}  
Then, we take $U^{n+1}=V^{n+1}$. The nonlinear subsystem (\ref{non_linear_ADI}) can be analytically integrated \cite{sheik_zhao_2020} so that the nonlinear instability is avoided. The linear subsystem (\ref{diffusion_ADI}) subject to the same jump and boundary conditions (\ref{eq:TRPBE_jump1}) - (\ref{eq:bc}) is first integrated in time by the implicit Euler scheme. Then the linear system involving three-dimensional (3D) unknowns is decomposed into one-dimensional (1D) linear algebraic systems \cite{sheik_zhao_2020} 
\begin{eqnarray}
		\left(1- \Delta t \delta_x^2\right)v_{i,j,k}^{*}&=&\left[1+ \Delta t \left(\delta_y^2+\delta_z^2 \right)\right]v_{i,j,k}^{n},\label{GFM-ADI}\\
		\left(1-\Delta t \delta_y^2\right)v_{i,j,k}^{**}&=&v_{i,j,k}^{*}- \Delta t \delta_y^2\left(v_{i,j,k}^{n}\right),\label{GFM-ADI2}\\
		\left(1- \Delta t \delta_z^2\right)v_{i,j,k}^{n+1}&=&v_{i,j,k}^{**}- \Delta t \delta_z^2\left(v_{i,j,k}^{n}\right).\label{GFM-ADI3}
\end{eqnarray}
where $\delta_x^2$, $\delta_y^2$ and $\delta_z^2$ are the central finite difference operators for the $x$, $y$, and $z$ directions, respectively. For example, away from the interface $\Gamma$, we have
$\delta_x^2 (v_{i,j,k}^n) = \frac{\epsilon_{i,j,k}}{h^2} (v_{i-1,j,k}^n-2v_{i,j,k}^n+v_{i+1,j,k}^n)$.
Each ADI equation has a tridiagonal structure, and can be efficiently solved by the Thomas algorithm with a complexity like $O(N_x \log N_x)$. Thus, for each time step, the complexity of the entire ADI time splitting scheme is on the order of $O(N \log N)$, where $N=N_x \times N_y \times N_z$ is the spatial degree of freedom. 

In the LOD scheme, the pseudo-time PBE (\ref{eq:TRPBE}) is split into three stages
\begin{eqnarray}
  \frac{\partial w}{\partial t}&=& -\frac{1}{2}\kappa^2 \sinh(w) \text{ with } W^n=U^n\text{ and } t \in \left[t_n,t_{n+1}\right]\nonumber,\\
 \frac{\partial v}{\partial t}&=&  \nabla \cdot (\epsilon\nabla v) \text{    with } V^n=W^{n+1}\text{ and } t \in \left[t_n,t_{n+1}\right]	 \label{diffusion_LOD-CN},\\
 \frac{\partial \tilde{w}}{\partial t}&=& -\frac{1}{2}\kappa^2 \sinh( \tilde{w}) \text{ with } \tilde{W}^n=V^{n+1}\text{ and } t \in \left[t_n,t_{n+1}\right].\nonumber
\end{eqnarray}
Then we have $U^{n+1}=\tilde{W}^{n+1}$. The two nonlinear subsystems are analytically integrated, while the linear subsystem subject to (\ref{eq:TRPBE_jump1}) - (\ref{eq:bc}) is decomposed into 1D diffusion equations, which are then discretized by the Crank-Nicolson scheme \cite{sheik_zhao_2020}
\begin{eqnarray}
	\begin{aligned}
	\left(1-\frac{\Delta t}{2}\delta_x^2\right)v^*_{i,j,k}&=\left(1+\frac{\Delta t}{2}\delta_x^2\right)v^n_{i,j,k},\\
	\left(1-\frac{\Delta t}{2}\delta_y^2\right)v^{**}_{i,j,k}&=\left(1+\frac{\Delta t}{2}\delta_y^2\right)v^*_{i,j,k},\\
	\left(1-\frac{\Delta t}{2}\delta_z^2\right)v^{n+1}_{i,j,k}&=\left(1+\frac{\Delta t}{2}\delta_z^2\right)v^{**}_{i,j,k}.\label{LODCN_eq}
	\end{aligned}
\end{eqnarray}
This method is called as the LODCN in \cite{sheik_zhao_2020}, and will be simply referred to as the LOD scheme in the present study. 

\subsection{Modified ghost-fluid method}
Because the potential $u$ and its flux are discontinuous across the interface $\Gamma$, jump conditions (\ref{eq:TRPBE_jump1}) and (\ref{eq:TRPBE_jump2}) have to be satisfied numerically in order to guarantee spatial convergence. Moreover, in the ADI and LOD time discretizations, the 3D problem is decomposed into 1D ones in order to achieve high efficiency. Most existing interface treatments in the literature cannot be applied to the present problem, because jump conditions have to be enforced in a tensor product decomposition manner to fit with the ADI/LOD framework. There exists a few matched ADI methods \cite{Zhao15,Li20} that deliver a second order of spatial accuracy for solving the parabolic interface problems in the ADI framework, which, however, break the symmetry of the finite difference matrix.  The greatest contribution of Ref. \cite{sheik_zhao_2020} is the introduction of a modified ghost-fluid method (GFM), which not only handles the discontinuous jumps in discretization, but also maintains a symmetric and tridiagonal matrix structure for each ADI/LOD step. The key challenge here lies in the flux jump condition (\ref{eq:TRPBE_jump2}), i.e., the normal direction is not a Cartesian direction so that its direct discretization naturally couples all Cartesian directions. In \cite{sheik_zhao_2020}, this difficulty is overcame through an approximation of (\ref{eq:TRPBE_jump2}) by three tensor product decomposed jump conditions in Cartesian directions 
\begin{equation}\label{GFM-jump}
\left[\epsilon \frac{\partial u}{\partial x}\right]\approx\epsilon^- \frac{\partial G}{\partial x}, \quad 
\left[\epsilon \frac{\partial u}{\partial y}\right]\approx\epsilon^- \frac{\partial G}{\partial y}, \quad \text{ and } \quad \left[\epsilon \frac{\partial u}{\partial z}\right]\approx\epsilon^- \frac{\partial G}{\partial z}.	  	
\end{equation}
In this manner, for each intersection point between the interface $\Gamma$ and a Cartesian grid line (say $x$-line), two 1D jump conditions are available, e.g. $[u]=G$ and $[\epsilon \frac{\partial u}{\partial x} ] = \epsilon^- \frac{\partial G}{\partial x}$. These two 1D conditions can be enforced to modify the finite difference operator $\delta_x^2$ near the interface, without breaking the symmetry \cite{sheik_zhao_2020}. The modified GFM yields a first order global accuracy in space. 

\subsection{Pseudo-time GFM algorithms}
By combining the ADI/LOD time integration with the GFM spatial discretization, two pseudo-time GFM algorithms will be investigated in this paper, i.e., the GFM-ADI and GFM-LOD schemes  \cite{sheik_zhao_2020}. The PBE solvers are usually benchmarked by the  electrostatic free energy $E_{\rm sol}$, which is defined as the energy released when the solute in free space is dissolved in the solvent. In the regularization methods \cite{Geng17}, one does not calculate the energy for the solute in the vacuum state. Instead, one just solves one regularized PBE in the water state, and $E_{\rm sol}$ can be approximated as
\begin{equation}\label{eq:solvation}
	E_{\rm sol}= \frac{1}{2}k_B T \sum\limits_{i=1}^{N_c}q_i\phi_{RF} ({\bf r}_i),
\end{equation}
where $\phi_{RF}$ is the reaction-field potential. In the pseudo-time GFM method \cite{sheik_zhao_2020}, the reaction-field potential $\phi_{RF}({\bf r})$ inside the solute domain $\Omega^-$ is simply the steady-state solution of $u({\bf r},t)$. 

In pseudo-time integration, either a zero solution or the solution of the linearized PBE is chosen as the initial solution $u({\bf r},0)$. Then one will solve the time-dependent PBE for a long time or until $t=T_{end}$ to ensure that the steady state is achieved. The steady-state convergence can also be checked by the energy difference. In particular, define $E^n_{\rm sol}$ to be the free energy at time step $t_n$. We calculate the energy difference at each time step
\begin{equation}\label{dE}
\Delta E^n_{\rm sol} = |E^n_{\rm sol} - E^{n-1}_{\rm sol}|.
\end{equation}
The convergence is assumed to be attained if the energy difference is less than a tolerance, i.e., $\Delta E^n_{\rm sol} < TOL$. Because the GFM discretization guarantees that each 1D finite difference matrix maintains symmetry, diagonally dominate, and tridiagonal \cite{sheik_zhao_2020}, the stability of the ADI and LOD methods is much better than those in the previous pseudo-time approaches \cite{Geng13,Zhao14,Leighton_Zhao}. The fully-implicit GFM-ADI method allows the use of a large time increment, which is efficient enough for most steady-state simulations, while the GFM-LOD method is unconditionally stable in time integration \cite{sheik_zhao_2020}.

Numerical experiments in \cite{sheik_zhao_2020} indicate that both GFM-ADI and GFM-LOD schemes are first order accurate in time. The GFM-ADI is more accurate when a smaller $\dt$ is used, while for large $\dt$ values, GFM-LOD is better. In terms of spatial accuracy, both methods achieve first order convergence in $L_\infty$ norm, and the orders in $L_2$ norm could be close to two. For each time step, the complexity of both methods scales like $O(N \log N)$ with $N$ being the number of total spatial unknowns. However, a constant $\dt$ is employed in \cite{sheik_zhao_2020}, so that many steps of time integration are commonly required to reach the steady-state. In \cite{sheik_zhao_2020}, the molecular surface is calculated based on the triangulation from the MSMS package, and is converted via a Lagrangian-to-Eulerian algorithm. Such a conversion is slightly simpler than that in the MIB method \cite{Yu2007}, because in the GFM, only interface locations are needed in (\ref{GFM-jump}), while the normal directions are not required. Nevertheless, the MSMS is known to induce instabilities in the PBE solution of some systems, specially when the density is large.

\section{Eulerian Solvent Excluded Surface (ESES)}
In this section, the implementation details of the ESES \cite{ESES,ESES_parallel} in the pseudo-time GFM approach \cite{sheik_zhao_2020} are reported. Numerical experiments are carried out to examine the stability and performance of the improved pseudo-time solver over the original one based on the MSMS. 

\subsection{ESES algorithm}
The Eulerian Solvent Excluded Surface (ESES) algorithm constructs the SES analytically over Cartesian grids \cite{ESES}. ESES takes as input the location of atom centers and radii, and calculates the analytical SES patches given by \cite{Connolly_surface}. Then ESES marks all points contained in the solvent accessible surface (SAS) as uncertain and those outside the SAS as outside the SES. ESES creates several auxiliary geometric features for the analytic patches. Membership or non-membership to these auxiliary features and the analytic patches are used to determine the classification of the majority of uncertain points. The remaining uncertain points are then classified by the number of times the SES intersects a grid line containing the point \cite{ESES}. This classification is an Eulerian representation of the SES. Finally, ESES computes distances between grid points and surface intersections as well as the normal direction of the surface at the intersection.

It has been observed that as MSMS surface density increases, the MSMS surface visually approaches that defined by ESES \cite{ESES}. On the other hand, for the electrostatic free energy calculated by the MIB-PBE algorithm \cite{Chen11}, the MSMS energy approaches that of the ESES as the density becomes larger. Moreover, the ESES is free of instability issues associated with the MSMS and its usage in finite difference methods avoids the need of the Lagrangian-to-Eulerian conversion \cite{Yu2007}. These factors motivate us to replace the MSMS by ESES in the pseudo-time GFM approach. 

\subsection{ESES implementation}
In this subsection we briefly describe the input and output of ESES and its integration with our source code. The input for ESES is the location of atom centers and radii. Our software requires the atom positions, radii, and charges, which we format in one file. We modified the ESES input method to accept the file that includes the atom charge, which is not used in surface calculation. This removes some duplicated data, as previously the radii and charges were stored separately, each with a copy of the coordinates.

ESES outputs a bounding box, the grid information, and the intersection information. The bounding box describes the number of points in each dimension, as well as the start and end points. The grid information classifies all points as inside or outside the molecule. The intersection information contains every intersection of the grid lines with the SES, which is described by the grid indices on either side of the surface, the distance from the inside grid point to the surface, the normal direction of the surface at the intersection, and the atom indices for the patch that was intersected.

Although there is a parallelized version of ESES \cite{ESES_parallel}, we use it as a single-threaded application. Utilizing the parallelized version would make surface computation much quicker for very large proteins. We select either the MSMS or ESES surface by command-line parameter. We set the ESES probe radius to be $r_p=1.4$\AA. ESES allows an extension parameter that describes the distance between the outermost atoms in each direction and the bounding box. We set this as the floor of $2r_p$, in order to approximate the bounding box we use for the MSMS surface. The use of an extension parameter is necessary to achieve the correct accuracy.

ESES is invoked via a system call from the Fortran code. In our tests of different pseudo-time algorithms, the same protein and spatial mesh are re-used. For such problems, re-computing the ESES surface for every test uses unnecessary computation time. To conserve time for our simulations that did not measure CPU time, we save the outputs in a directory categorized by PDB ID and grid size. We then read in the bounding box and grid information to set up our grid. In our previous MSMS package, we encoded points inside the MSMS surface as $-1$ and outside as $1$, while ESES does the opposite. Upon input, we negate the ESES values to avoid changing well-tested code.

Next, we read the intersection information. We use this to classify the irregular points, those with a neighboring grid point on the other side of the surface, which are used in the GFM described in \cite{sheik_zhao_2020}. We read the file once to get the number of intersections and estimate the maximum number of irregular points as no more than $3$ times the number of intersections. We re-read the intersection information into arrays declared for this maximum number. This method of input for irregular points is not efficient in CPU or memory usage, and is a point for future improvement.

\subsection{ESES validation}
In this subsection, we validate the ESES surface in the pseudo-time GFM package, and compare it with the original one based on the MSMS surface. Since the visual difference between the ESES and MSMS has been studied in \cite{ESES}, in the present study, we will focus on the difference in electrostatic free energy calculated by the ESES and MSMS, based on the same pseudo-time GFM algorithm. In all tests, the nonlinear PBE is solved with $\epsilon^+=80$ and $\epsilon^-=1$.  Two pseudo time methods, i.e., ADI and LODCN schemes \cite{sheik_zhao_2020} will be employed, with the latter being referred to as the LOD for simplicity. Throughout this paper, the length is reported with units of \AA~and the electrostatic free energy has the unit of kcal/mol. For all pseudo-time simulations in this section, the computation will stop if either the time $t$ reaches $T_{end}=10$ or the difference in free energy satisfies $E^n_{sol} < 10^{-4}$. We note that $T_{end}=10$ maybe not longer enough for the steady-state convergence in protein studies. Nevertheless, this shortened computation will not affect our comparison between the MSMS and ESES, and the conclusion to be drawn. 

\begin{figure}[!htbp]
 \centering
 \includegraphics[width=\textwidth]{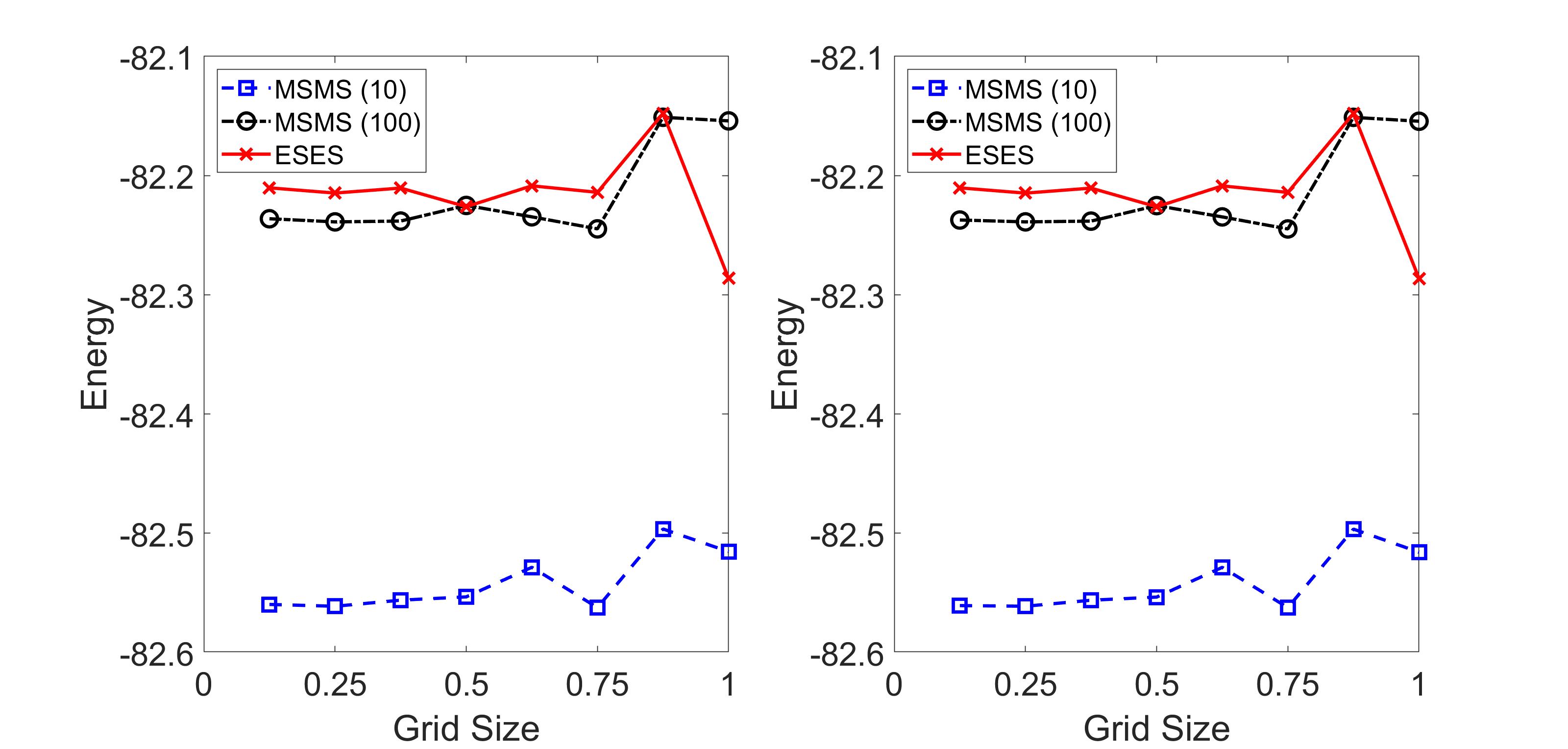}
 \caption{Convergence in electrostatic free energies for one atom model by using varying grid sizes based on an ESES surface, an MSMS surface with density 10, and an MSMS surface with density 100. Left: GFM-ADI; Right: GFM-LOD.}
 \label{fig:1ato_spatial}
\end{figure}
\subsubsection{Kirkwood sphere with one centered charge}
We consider a spherical cavity with radius $R=2$ and one centered charge $q=1$.  The ionic strength is chosen such that the Debye-Huckle constant $\kappa=1$. The analytical energy for the one atom model is available for both Poisson's equation and linearized PBE, while the analytical energy for the nonlinear PBE is unavailable, but is known to be very close to that of the Poisson's equation \cite{Geng17}, i.e., $\Delta G = -81.9782$ kcal/mol. In \cite{sheik_zhao_2020}, in order to fairly test the accuracy of the pseudo-time GFM algorithm, an exact sphere is used as the molecular surface. In the present study, in order to compare the ESES with MSMS, the molecular surface is generated numerically. Visually, we have verified that the ESES surface converges to the sphere as the grid size $h$ goes to zero. The default density of the MSMS surface is 10. By increasing it from 10 to 100, the MSMS surface becomes closer to the sphere.

\begin{table}[!htbp]
\centering
\begin{tabular}{|l|l|l|l|l|l|}\hline
Method & $h$     & ESES         & MSMS(10)      & MSMS(100)     & Sphere \cite{sheik_zhao_2020}      \\\hline
\multirow{5}{*}{GFM-ADI}  
& 1     & -82.28624893 & -82.51547246 & -82.1541084  & -82.132181 \\
& 0.5   & -82.22609329 & -82.55340833 & -82.22504644 & -82.063724 \\
& 0.25  & -82.2146617  & -82.56126599 & -82.23883839 & -82.051117 \\
& 0.125 & -82.21027973 & -82.55973867 & -82.23620811 & -82.046462 \\
\hline
\multirow{5}{*}{GFM-LOD}  
& 1     & -82.28627562 & -82.51588465 & -82.15442576 & -82.132148 \\
& 0.5   & -82.22615189 & -82.55359102 & -82.22520013 & -82.063684 \\
& 0.25  & -82.21470823 & -82.56125163 & -82.23882793 & -82.051064 \\
& 0.125 & -82.21030578 & -82.56070777 & -82.23726394 & -82.046402\\\hline
\end{tabular}
\caption{Convergence in electrostatic free energy for one atom model. For the MSMS, two densities values 10 and 100 are studied. In the last column, the molecular surface is taken as the exact sphere \cite{sheik_zhao_2020}.}
\label{tab:1ato_spatial}
\end{table}

We first test the energy convergence by considering various values of $h$, see Fig.  \ref{fig:1ato_spatial}. Here the GFM-ADI and GFM-LOD algorithms are employed with $\dt=0.001$ and $T_{end}=10$. It can be seen that the energies of the ESES and MSMS with density $100$ converge to the same place. The MSMS surface with density $10$ is slightly offset from the other two, but appears to converge just as quickly to its final solvation free energy. Table \ref{tab:1ato_spatial} shows a subset of these results and also includes the energies reported in \cite{sheik_zhao_2020}  by using the exact sphere as the molecular surface. With $h \le 1$, the energies become pretty close to that of exact sphere. With the smallest $h=0.125$, the ESES and MSMS energy approaches to $-82.21$ and $-82.23$, respectively. These values are reasonably close to the limiting value of the exact sphere, i.e., $-82.05$. Nevertheless, we note that the difference between the limiting energies of ESES/MSMS and exact sphere does not seem to approach to zero as $h \to 0$. This perhaps suggests that the numerical molecular surface will always invoke certain approximation error in practice, even though it is negligibly small. 

\begin{figure}[!htbp]
 \centering
 \begin{subfigure}[b]{0.45\textwidth}
 \centering
 \includegraphics[width=\textwidth]{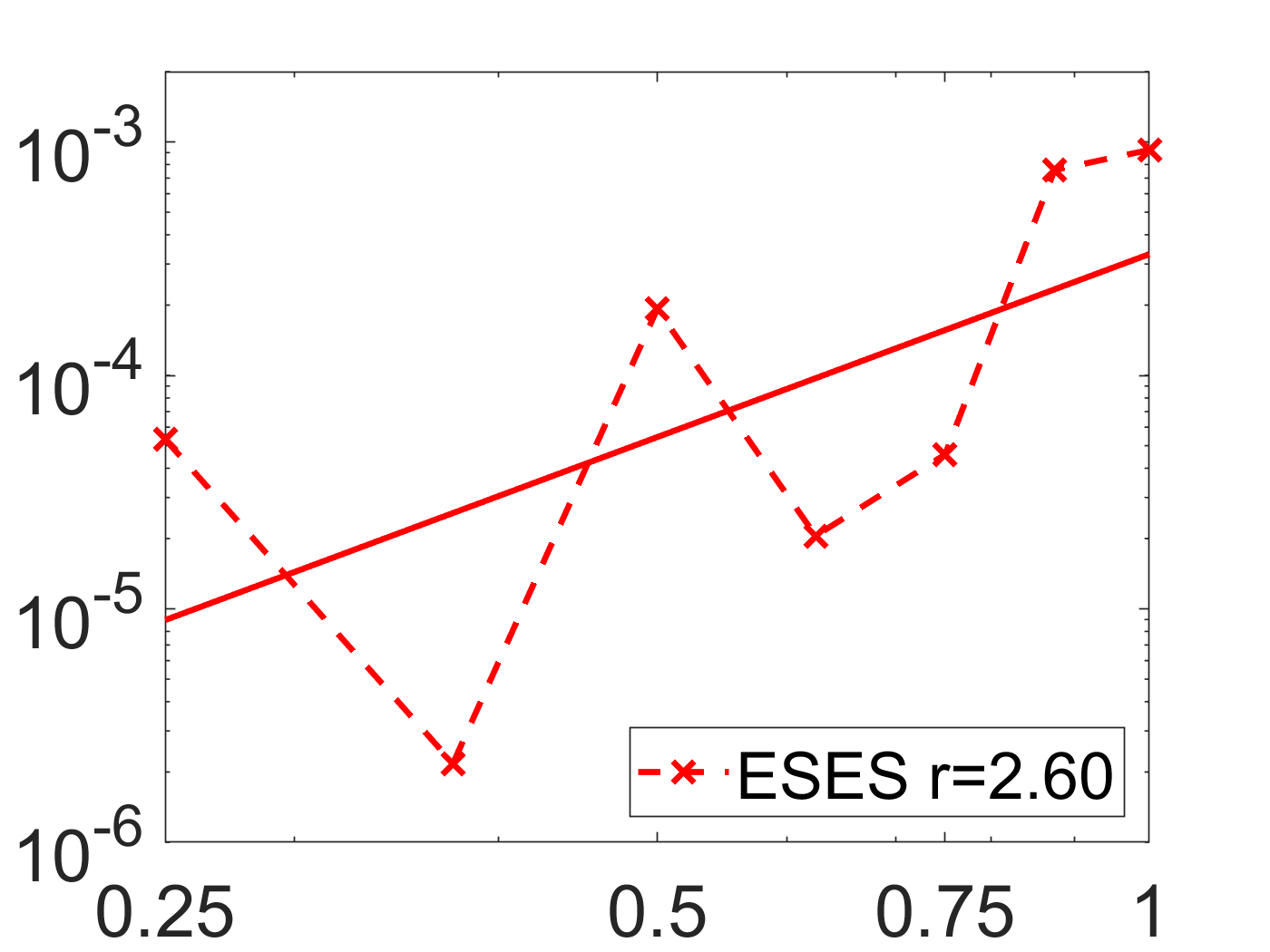}
 \caption{ESES GFM-ADI}
 \label{subfig:1ato_spatial_eses_gfmadi}
 \end{subfigure}
 \begin{subfigure}[b]{0.45\textwidth}
 \centering
 \includegraphics[width=\textwidth]{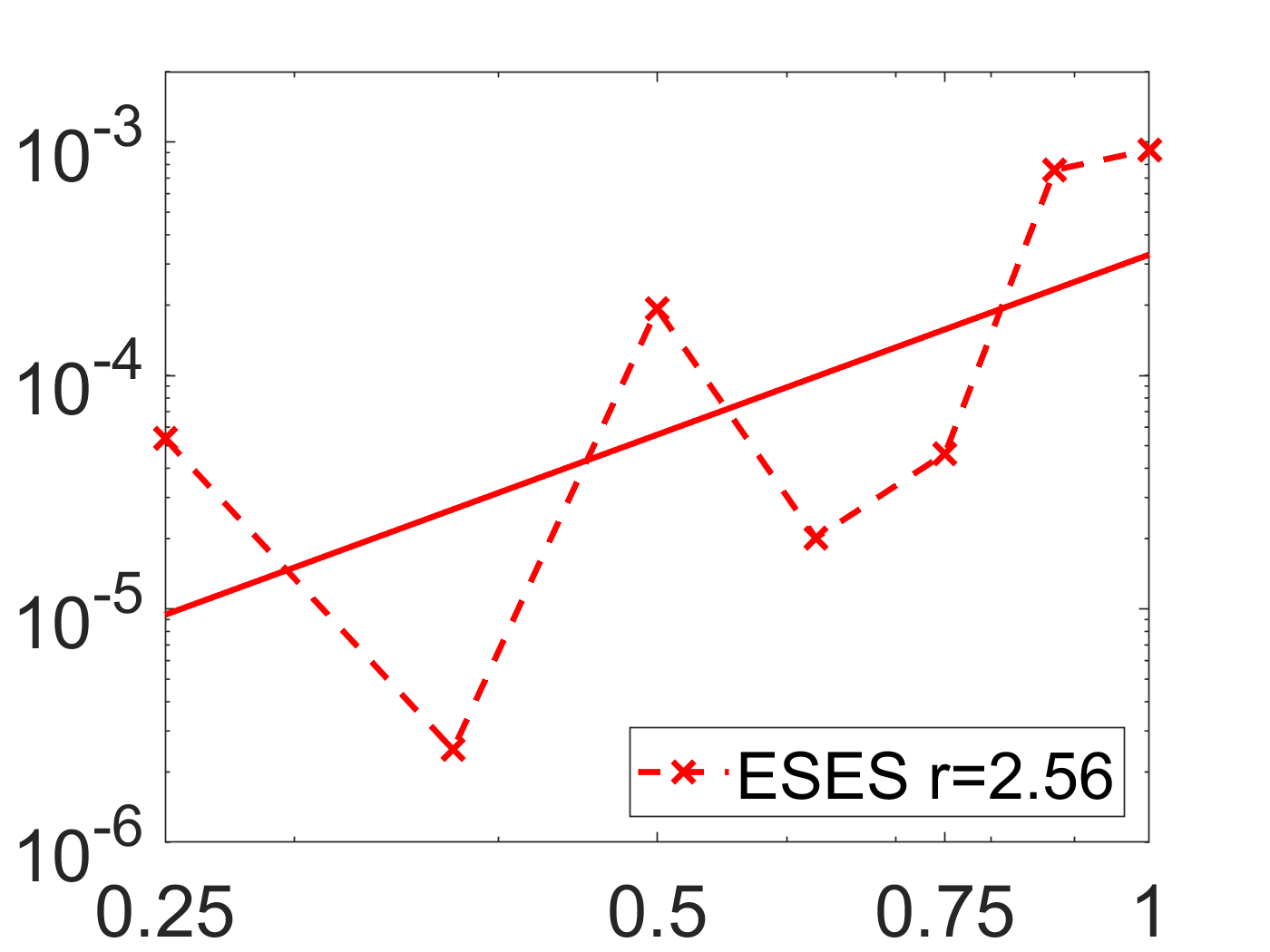}
 \caption{ESES GFM-LOD}
 \label{subfig:1ato_spatial_eses_gfmlod}
 \end{subfigure}
 \begin{subfigure}[b]{0.45\textwidth}
 \centering
 \includegraphics[width=\textwidth]{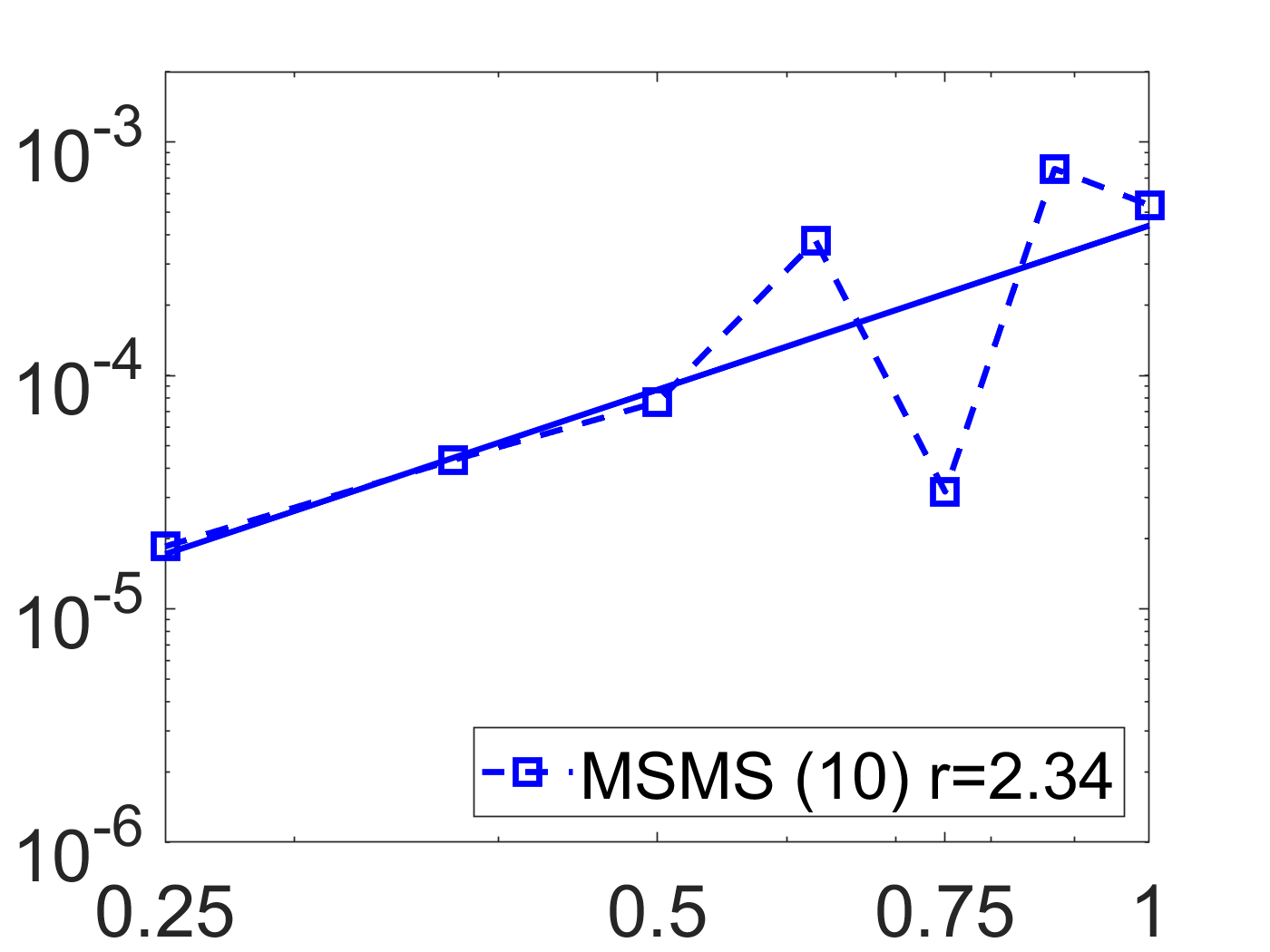}
 \caption{MSMS 10 GFM-ADI}
 \label{subfig:1ato_spatial_msms010_gfmadi}
 \end{subfigure}
 \begin{subfigure}[b]{0.45\textwidth}
 \centering
 \includegraphics[width=\textwidth]{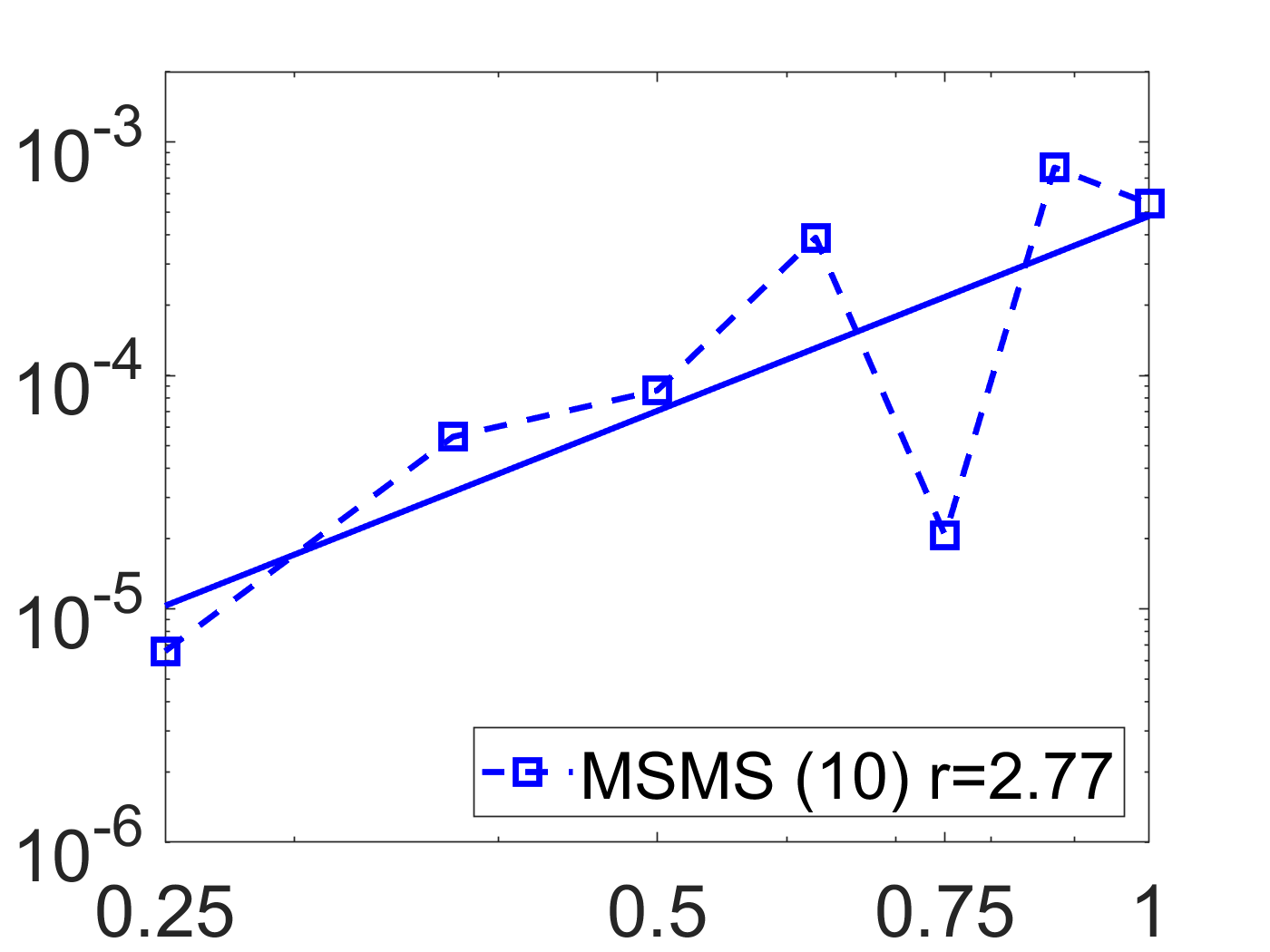}
 \caption{MSMS 10 GFM-LOD}
 \label{subfig:1ato_spatial_msms010_gfmlod}
 \end{subfigure}
 \begin{subfigure}[b]{0.45\textwidth}
 \centering
 \includegraphics[width=\textwidth]{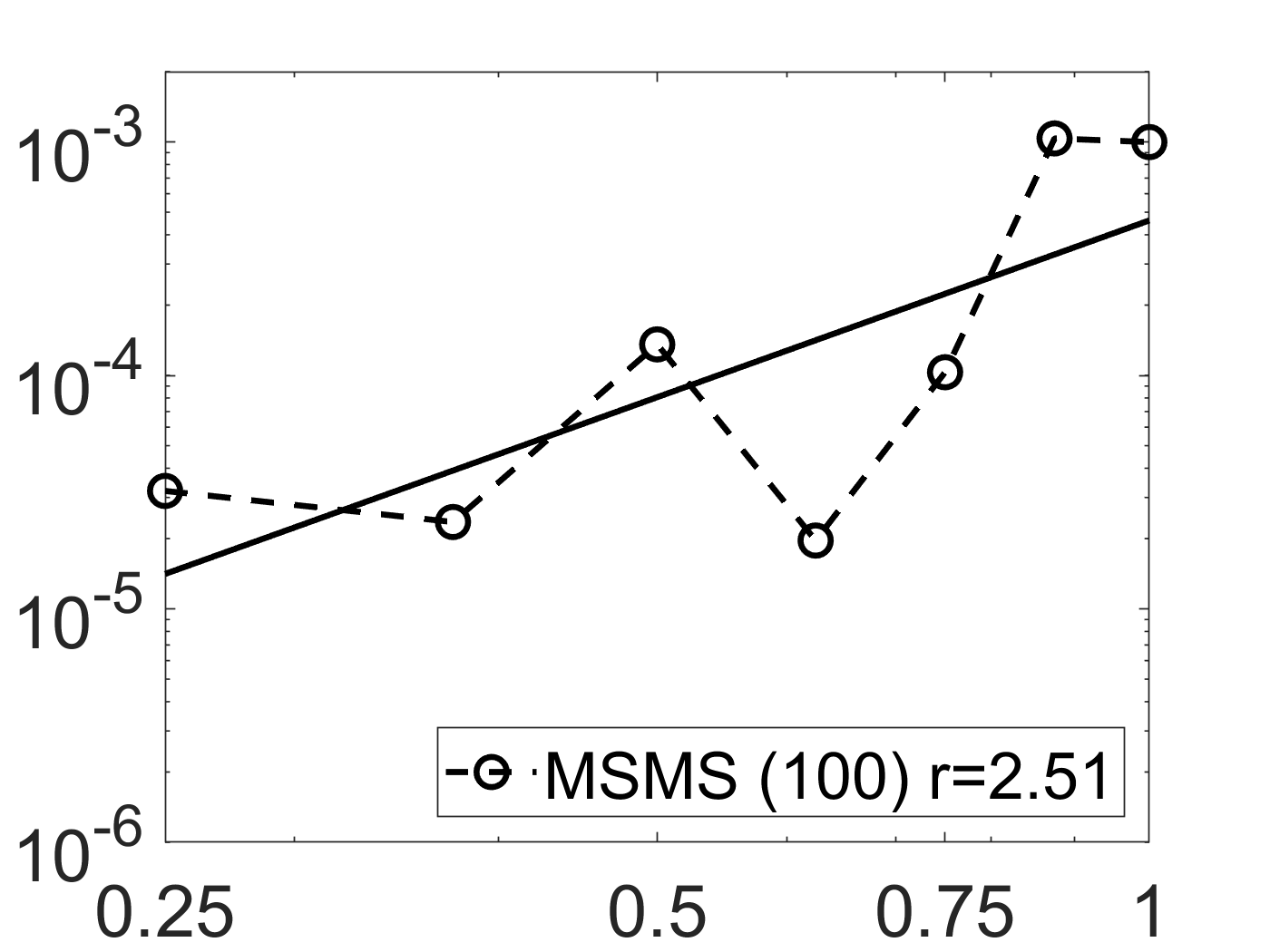}
 \caption{MSMS 100 GFM-ADI}
 \label{subfig:1ato_spatial_msms100_gfmadi}
 \end{subfigure}
 \begin{subfigure}[b]{0.45\textwidth}
 \centering
 \includegraphics[width=\textwidth]{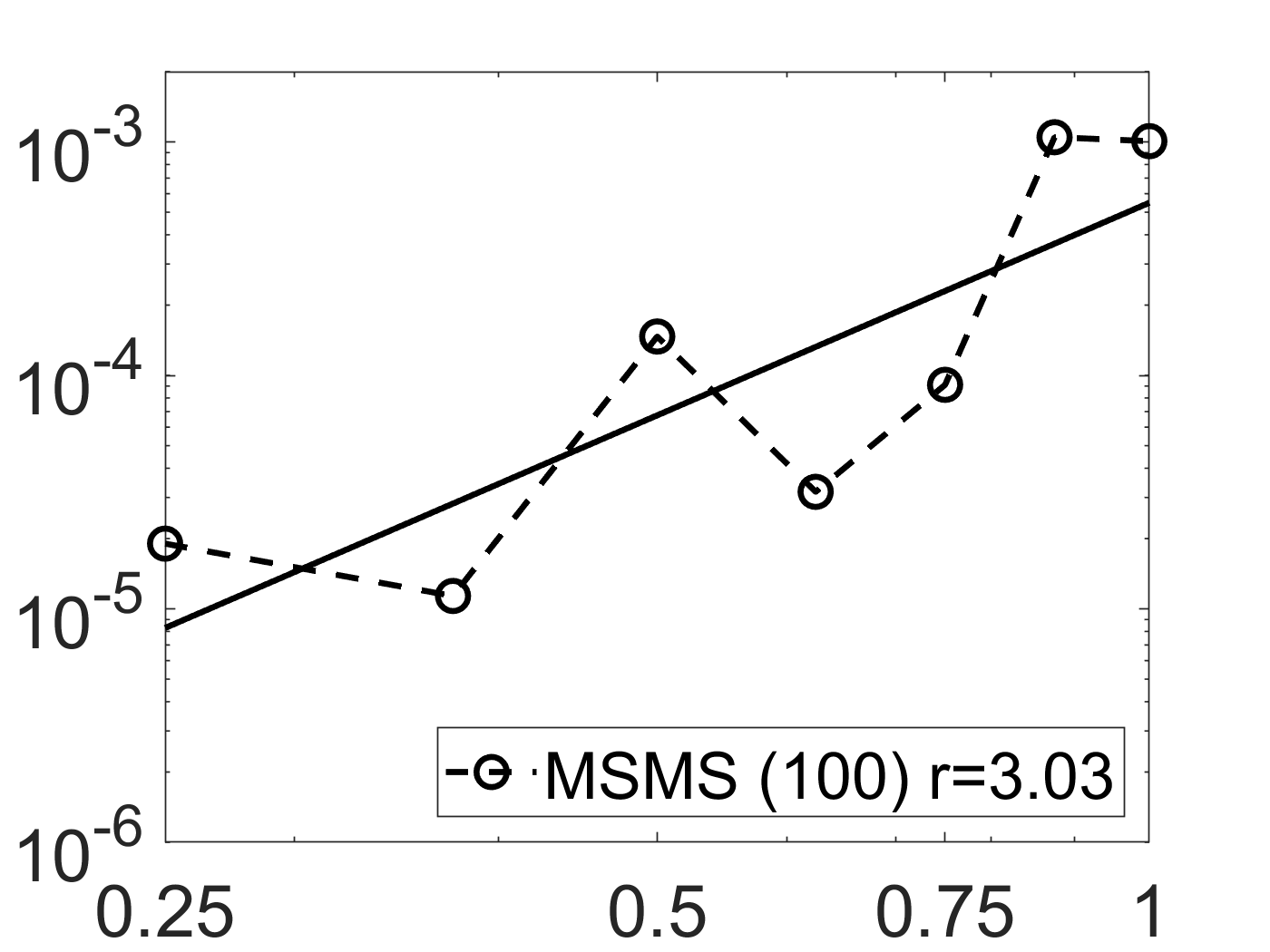}
 \caption{MSMS 100 GFM-LOD}
 \label{subfig:1ato_spatial_msms100_gfmlod}
 \end{subfigure}
 \caption{self-convergence study for the one atom model by taking the energy computed with $h=0.125$ as the reference solution for each case. Both GFM-ADI and GFM-LOD methods are studied with an ESES surface, an MSMS surface with density 10, and an MSMS surface with density 100. Also plotted are the least squares regression lines of best fit for the log of each relative error. Reported for each regression is the rate of convergence, $r$.}
 \label{fig:1ato_spatial_errors}
\end{figure}

In Fig. \ref{fig:1ato_spatial_errors}, self-convergence studies are carried out for both GFM-ADI and GFM-LOD methods based on the ESES, MSMS with density 10, and MSMS with density 100. For each case, we choose the reference energy to be the one obtained by $h=0.125$, and calculate the relative errors of other grid sizes with respect to the reference energy. By plotting the error against the grid size, we are able to numerically analyze the self-convergence order of the pseudo-time GFM algorithm. Such convergence rates are also shown in the legends of each subfigure, and they are all above two. We note that the self-convergence rate is not equivalent to the actual numerical order - the latter can be numerically detected only when analytical energy is available. But the self-convergence rate usually can imply the actual numerical order, especially when the grid size $h$ is sufficiently small. This study demonstrates that the GFM algorithm could achieve a spatially second order convergence in free energy calculation.

\begin{figure}[!htbp]
 \centering
 \begin{subfigure}[b]{0.45\textwidth}
 \centering
 \includegraphics[width=\textwidth]{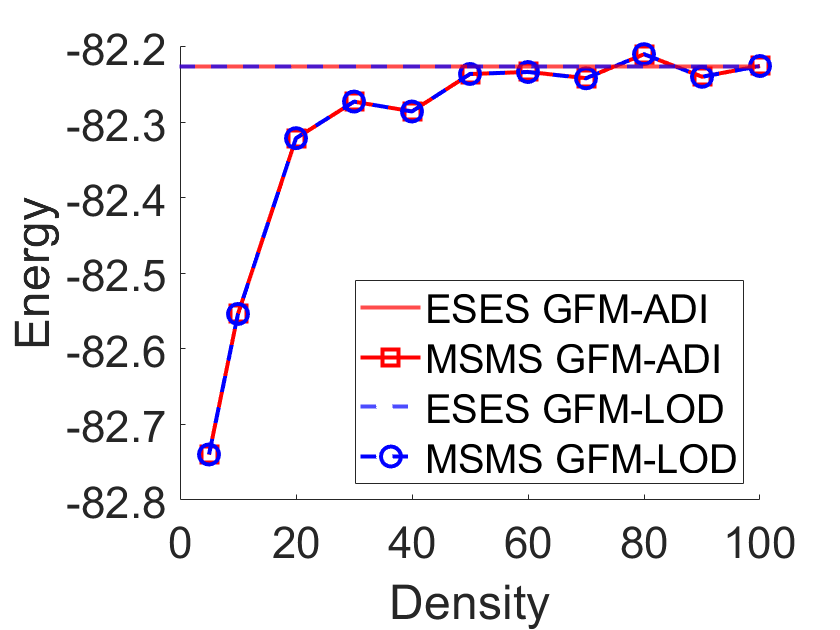}
 \caption{1 atom}
 \label{subfig:1ato_density}
 \end{subfigure}
 \begin{subfigure}[b]{0.45\textwidth}
 \centering
 \includegraphics[width=\textwidth]{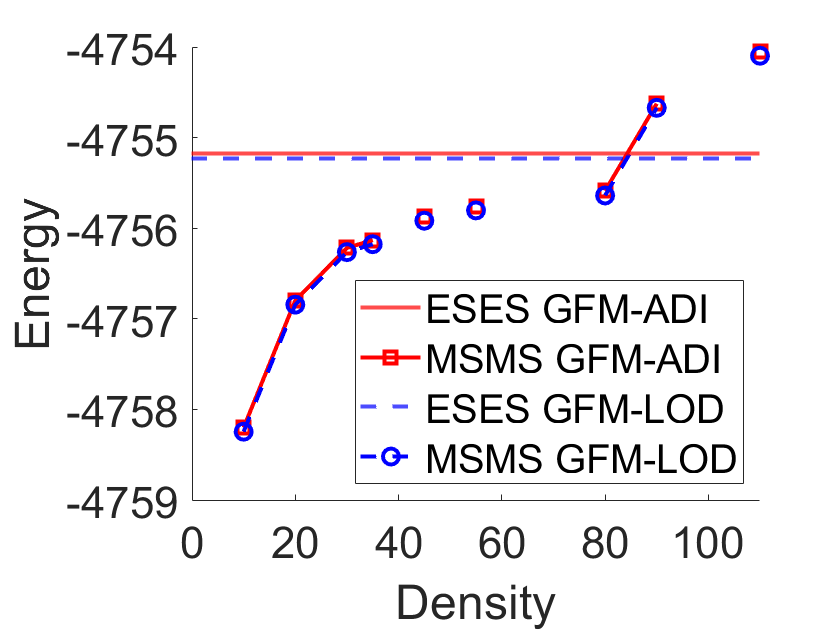}
 \caption{1a2e}
 \label{subfig:1a2e_density}
 \end{subfigure}
  \caption{The electrostatic free energy of the ESES and MSMS with different densities for (a) one atom model with $h=0.5$ and (b) protein 1a2e with $h=0.5$. In (b), both GFM solvers diverged for the MSMS surface of 1a2e with densities $40, 50, 60, 65, 70, 75, 100,$ and $105$. For density $95$, the converged value of both solvers involves an extremely high error, and is not shown in the figure. By increasing the surface density in both cases from $5$ to $100$, the MSMS energy approaches that of the ESES.}
 \label{fig:1ato_and_1a2e_density_convergence}
\end{figure}

We next demonstrate the convergence of the MSMS surface with respect to the density. By using $h=0.5$, the electrostatic free energy produced by different MSMS densities are depicted in Fig. \ref{fig:1ato_and_1a2e_density_convergence} (a). Obviously, the MSMS energy converges to the ESES energy as density becomes larger. When density is above 50, the MSMS energy oscillates around the ESES energy. 
This validates the accuracy of the ESES surface for the one atom model. We also note that by using a large density, the MSMS produces more triangles to represent the SES, which becomes numerically more demanding.

\begin{figure}[!htbp]
 \centering
 \includegraphics[width=\textwidth]{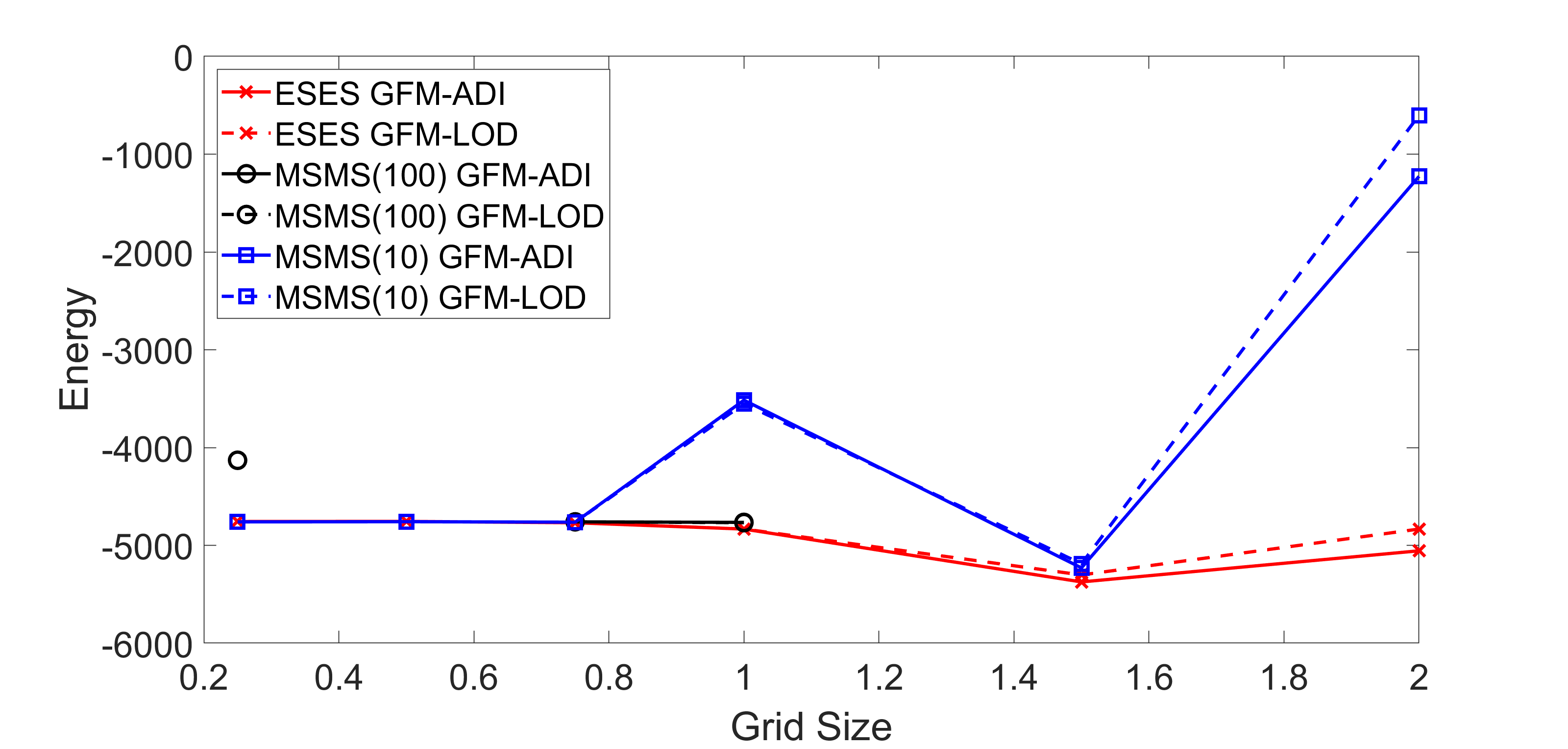}
 \caption{Electrostatic free energies of 1a2e when using varying grid sizes based on an ESES surface, an MSMS surface with density 10, and an MSMS surface with density 100. For the MSMS with density 100, the energies were non-negative for $h=1.5$ and $h=2.0$, while for $h=0.5$, both pseudo-time solvers diverged. These divergent values are not depicted in the figure.}
 \label{fig:1a2e_spatial}
\end{figure}

\begin{table}[!htbp]
\centering
\begin{tabular}{|l|l|l|l|l|}\hline
Method & $h$   & ESES     & MSMS(10)   &MSMS(100)     \\\hline
\multirow{4}{*}{GFM-ADI} & 1 & -4833.64 & -3515.85 & -4765.34    \\
& 3/4 & -4769.28 & -4763.13 & -4759.73    \\
 &1/2 & -4755.23 & -4758.26 & Fail         \\
& 1/4   & -4756.45 & -4760.05 & -4129.98    \\\hline
\multirow{4}{*}{GFM-LOD}& 1 & -4832.90 & -3552.67 & -4768.46    \\
& 3/4 & -4769.18 & -4763.17 & -4759.78    \\
& 1/2 & -4755.28 & -4758.30 & Fail         \\
& 1/4   & -4756.49 & -4760.08 & -4130.61    \\\hline
\end{tabular}
\caption{Convergence in electrostatic free energy for a protein (PDB ID: 1a2e).  For the MSMS, two densities values 10 and 100 are studied. }
\label{1a2e}
\end{table}

\subsubsection{Electrostatics on proteins}
We next test the pseudo-time GFM methods and ESES/MSMS molecular surface by studying several protein systems with complex geometry, and the protein structures were processed by using the CHARMM force field. The nonlinear PBE is solved with the ionic strength $I_s=0.15$ M. For both GFM-ADI and GFM-LOD schemes, we take $\dt=0.001$ and $T_{end}=10$. 

We first focus on a particular protein with PDB ID: 1a2e \cite{1A2E}. The electrostatic free energies calculated by the ESES, MSMS with density 10, and MSMS with density 100, are shown in Fig. \ref{fig:1a2e_spatial} for different grid size $h$. By changing $h$ from 2 to 0.25, the ESES energies for both GFM-ADI and GFM-LOD behave similarly and show a fast convergence. For the MSMS with the default density 10, the energies for a large $h$ obviously involve a quite large error. When $h<1$, the energies become visually the same as those of the ESES. For the MSMS with density 100, a lot of instability cases have been experienced. In particular, the energies were non-negative for $h=1.5$ and $h=2.0$, while for $h=0.5$, both pseudo-time solvers diverged. These divergent values are not depicted in Fig. \ref{fig:1a2e_spatial}. Other energies of the MSMS with density 100 are still problematic. For example, such values at $h=0.25$ obviously do not agree with the others. These instability issues are simply due to the use of a dense MSMS triangulation, which not only renders the molecular surface generation more complicated, but also brings more difficulty to the Lagrangian-to-Eulerian conversion. 

The electrostatic free energies calculated by three SES definitions are listed in Table \ref{1a2e} for $h \le 1$. For the ESES surface, as $h$ reduces from 1 to $0.25$, both the GFM-ADI and GFM-LOD results converge to the same place, which is about $-4756.5$. A fast self-convergence can be seen. For the MSMS with density 10, large errors are presented for $h=1$, but for smaller $h$ values, they rapidly converge to $-4760.1$. For the MSMS with density 100, the energies at $h=1$ and 0.75 seem to be on the right track for a convergence, but both GFM methods diverged at $h=0.5$ and the energies at $h=0.25$ are obviously wrong. Again, the MSMS energy is ruined by the instability introduced by the large density.

We next study the convergence of the MSMS surface with respect to the density. By using $h=0.5$, the electrostatic free energies produced by different MSMS densities are depicted in Fig. \ref{fig:1ato_and_1a2e_density_convergence} (b). For both GFM-ADI and GFM-LOD schemes, the MSMS energy shows a convergence trend initially. But the limiting values are hard to determine, because both solvers are unstable for so many large densities. This example demonstrates that for real proteins, the MSMS energy with a high density is not reliable, while the energy with a low density is not accurate enough. 

\begin{figure}[!htbp]
\centering
    \begin{subfigure}[b]{.75\textwidth}
        \includegraphics[width=\textwidth]{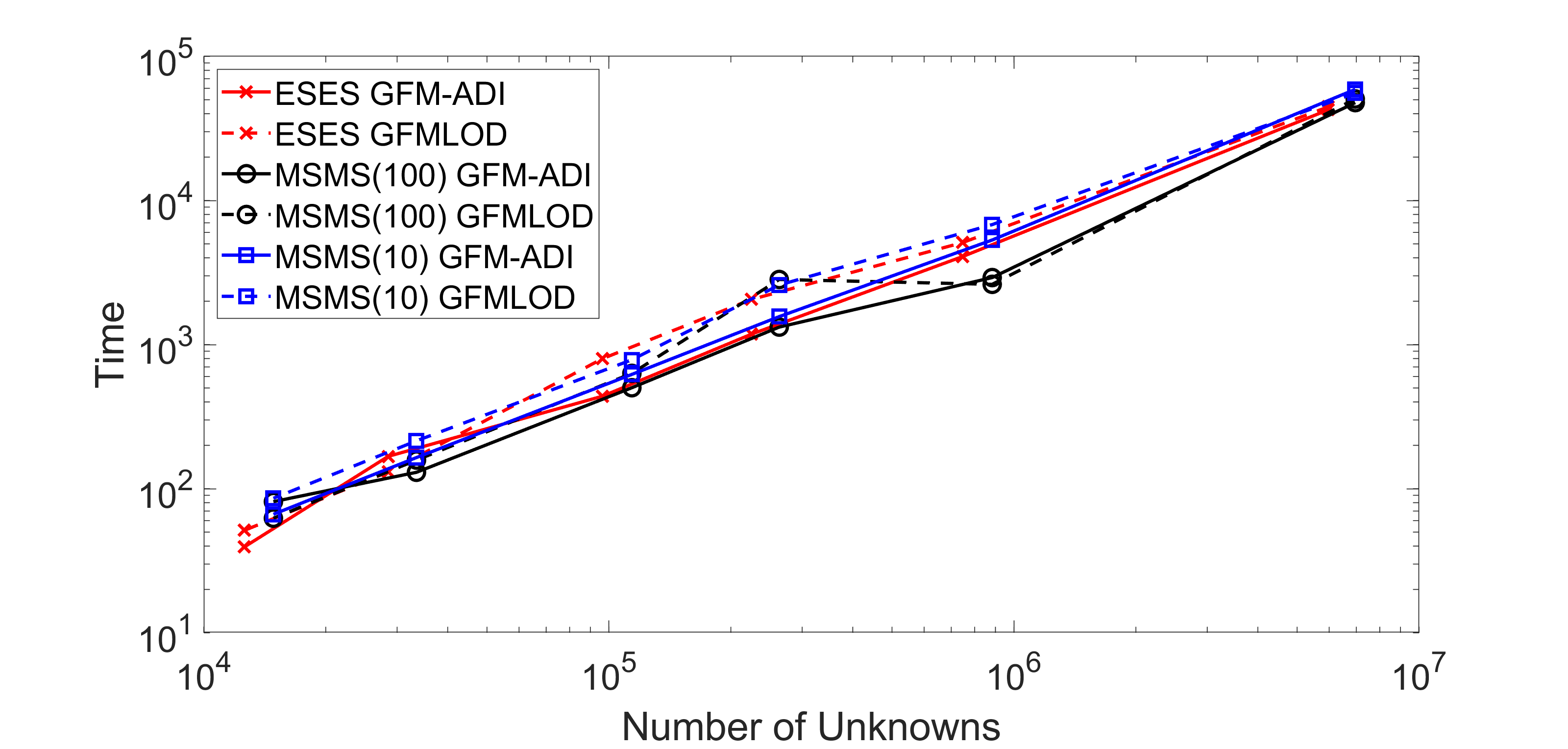}
        \caption{Time to solve the PBE}
        \label{fig:1a2e_nlpb_times}
    \end{subfigure}
    \begin{subfigure}[b]{.75\textwidth}
        \centering
        \includegraphics[width=\textwidth]{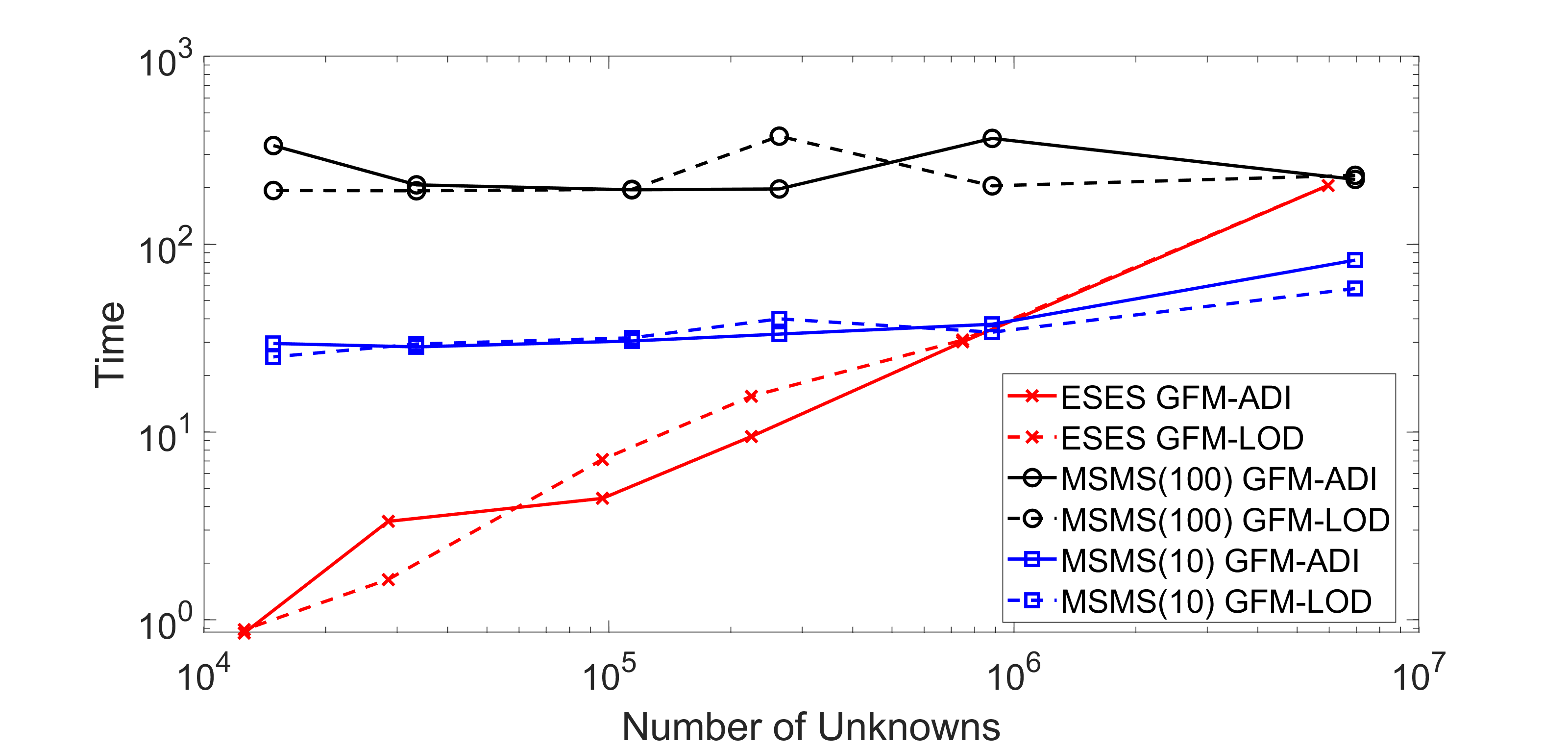}
        \caption{Time to compute the molecular surface}
        \label{fig:1a2e_surface_times}
    \end{subfigure}
    \begin{subfigure}[b]{.75\textwidth}
        \centering
        \includegraphics[width=\textwidth]{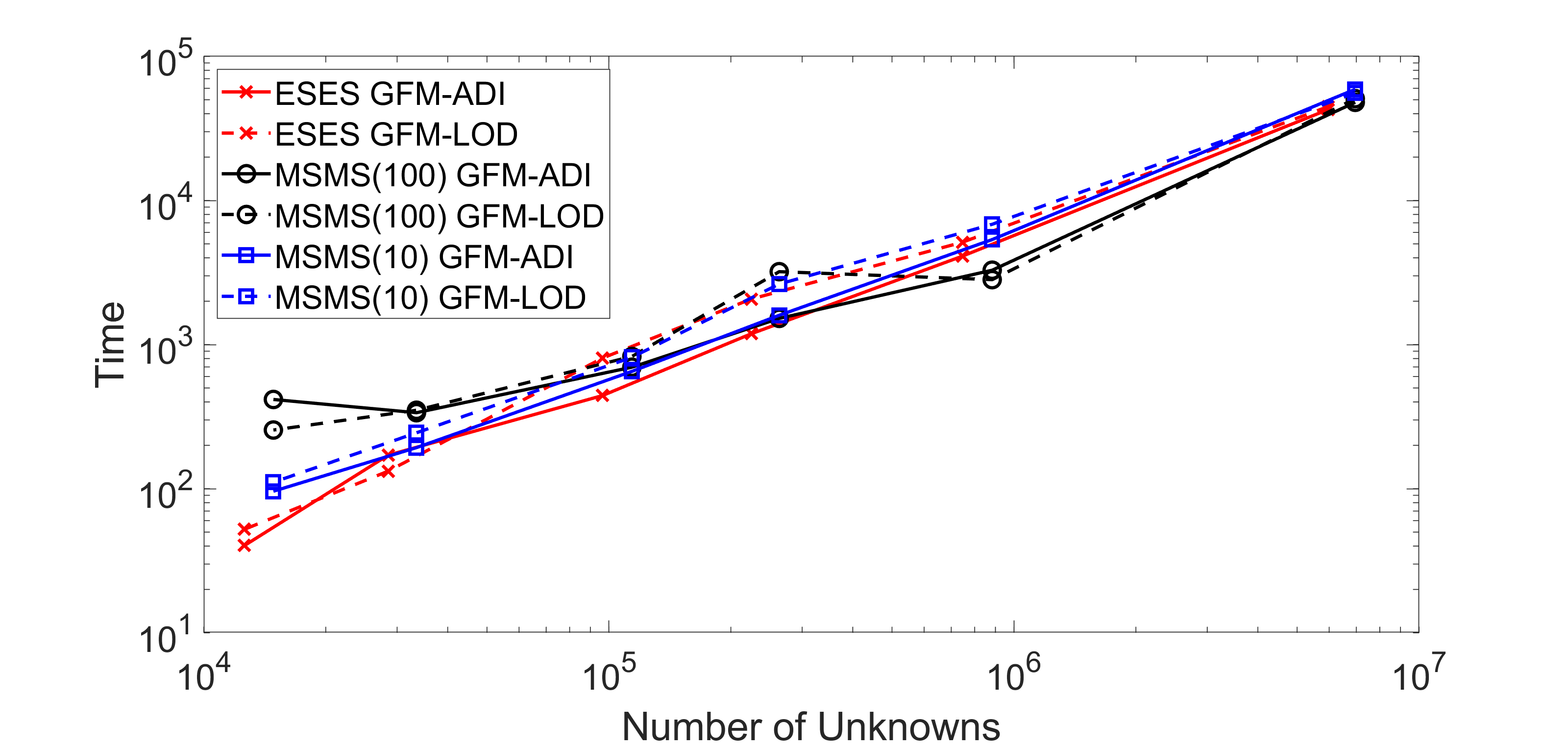}
        \caption{Total time to compute the molecular surface and solve the PBE.}
        \label{fig:1a2e_total_times}
    \end{subfigure}
\caption{CPU time in seconds taken for calculating electrostatic free energy of protein 1a2e by using varying grid sizes. Two pseudo-time solvers with a fixed $\dt=0.001$ and $T_{end}=10$ are considered and three surfaces, i.e., an ESES surface, an MSMS surface with density 10, and an MSMS surface with density 100, are employed. 
(a) is the time taken to solve the PBE alone, (b) is the time taken to compute the molecular surface alone, and (c) is the total time to compute the energy, including computing the molecular surface and solving the PBE.} 
 \label{fig:1a2e_spatial_times}
\end{figure}

\begin{figure}[!htbp]
\centering
    \begin{subfigure}[b]{0.45\textwidth}
        \centering
        \includegraphics[width=\textwidth]{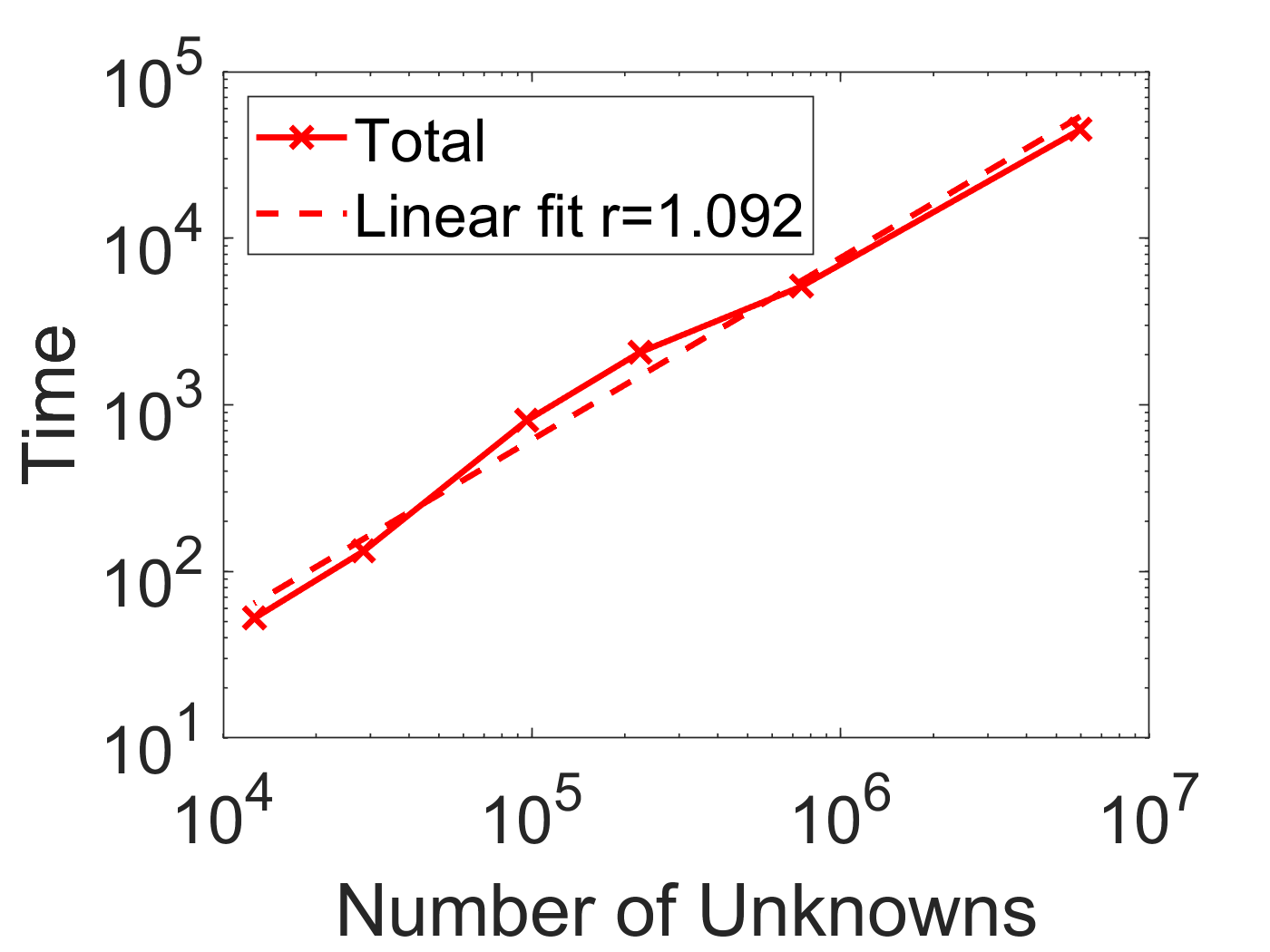}
        \caption{ESES GFM-ADI}
        \label{subfig:1a2e_spatial_time_gfmADI_eses}
    \end{subfigure}
    \begin{subfigure}[b]{0.45\textwidth}
        \centering
        \includegraphics[width=\textwidth]{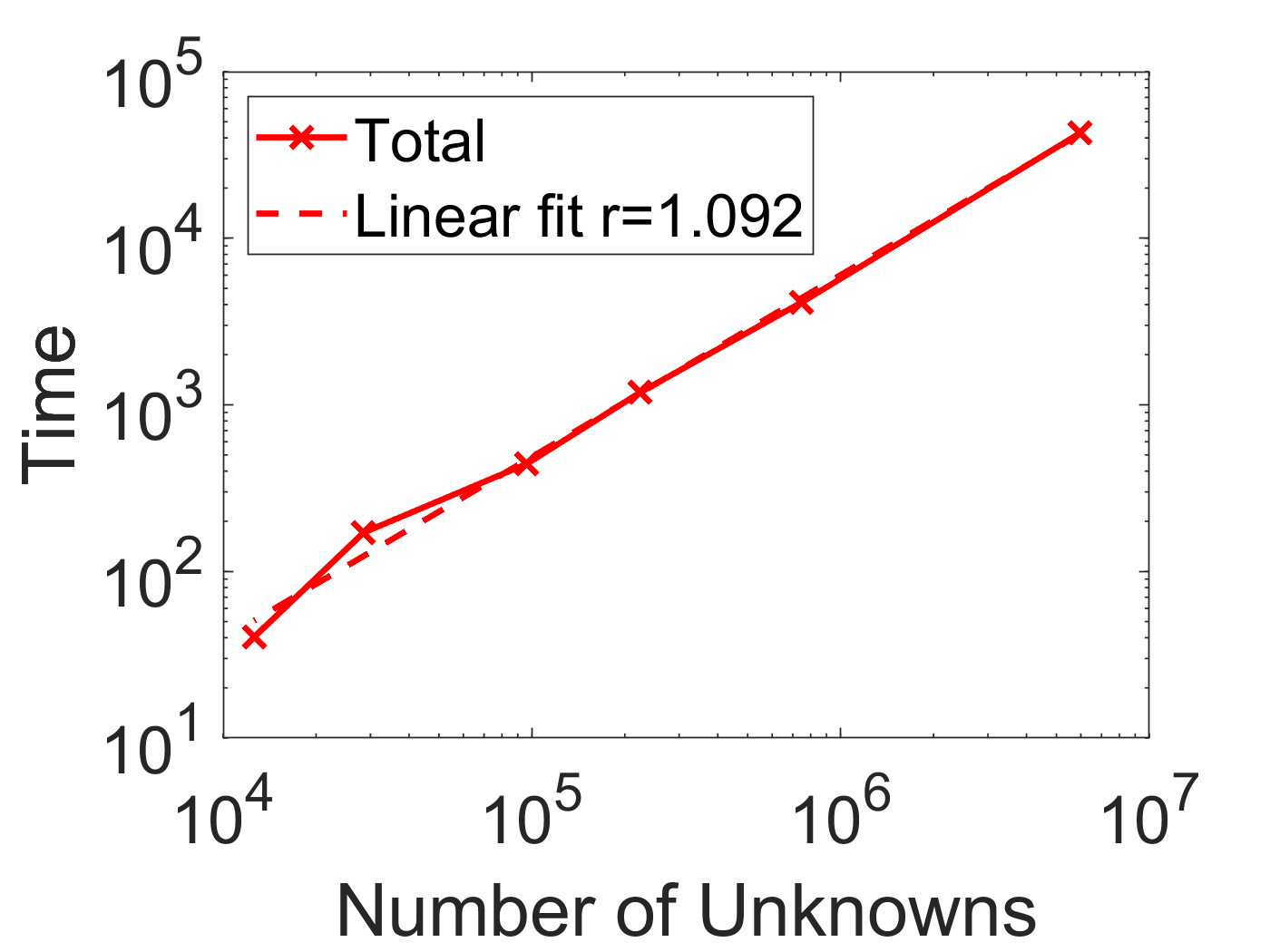}
        \caption{ESES GFM-LOD}
        \label{subfig:1a2e_spatial_time_gfmLOD_eses}
    \end{subfigure}
    \begin{subfigure}[b]{0.45\textwidth}
        \centering
        \includegraphics[width=\textwidth]{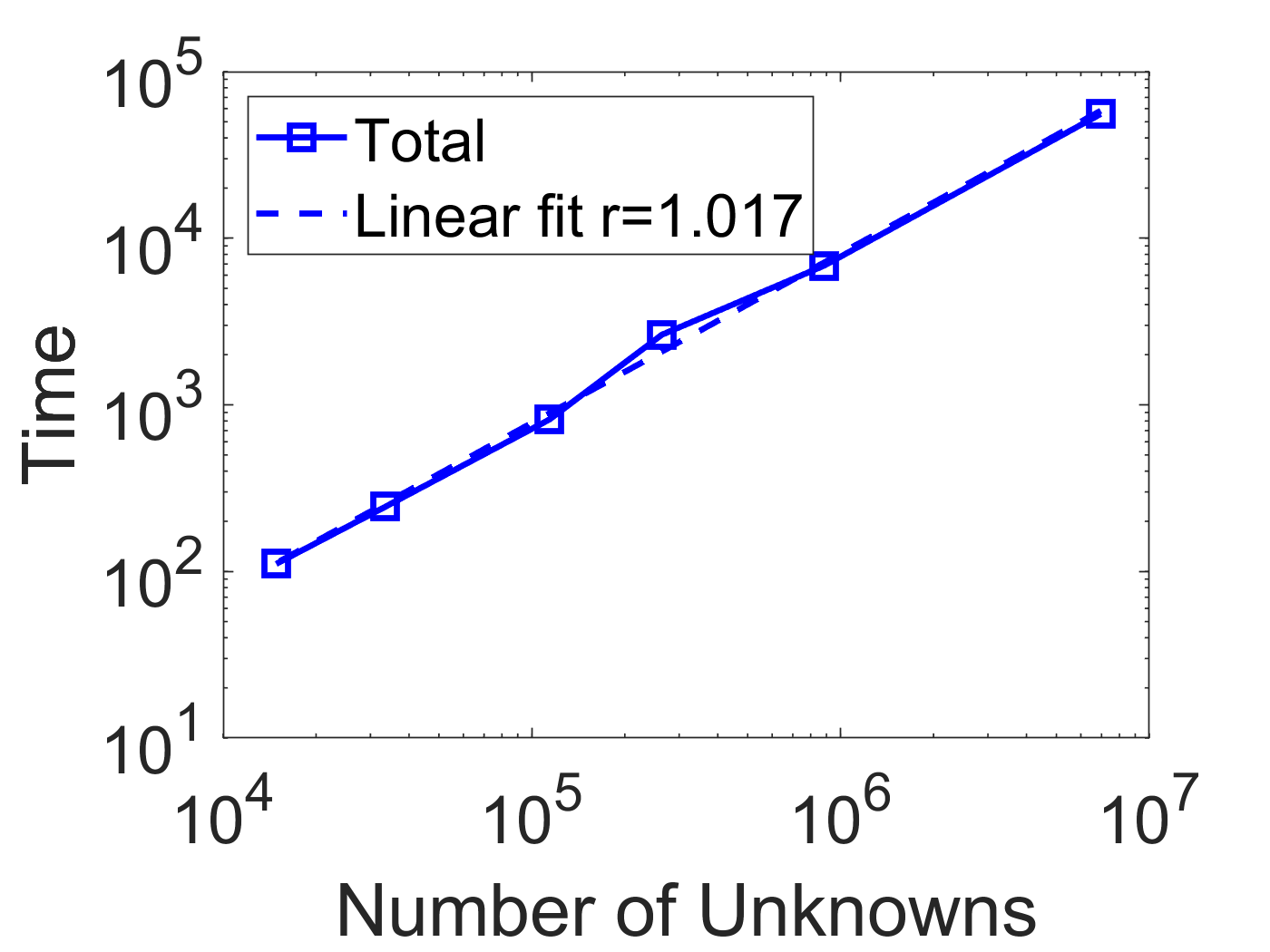}
        \caption{MSMS (10) GFM-ADI}
        \label{subfig:1a2e_spatial_time_gfmADI_msms10}
    \end{subfigure}
    \begin{subfigure}[b]{0.45\textwidth}
        \centering
        \includegraphics[width=\textwidth]{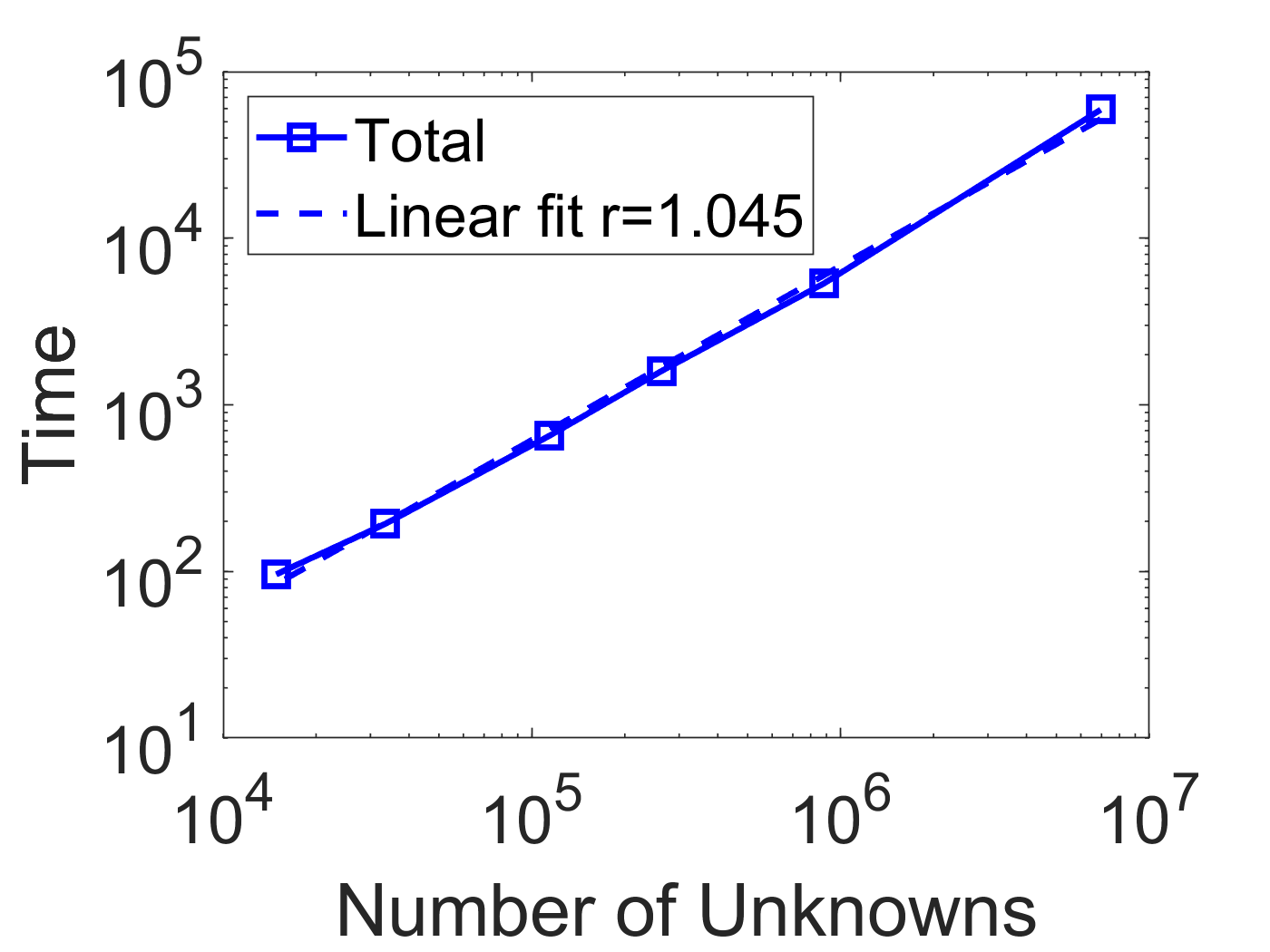}
        \caption{MSMS (10) GFM-LOD}
        \label{subfig:1a2e_spatial_time_gfmLOD_msms10}
    \end{subfigure}
    \begin{subfigure}[b]{0.45\textwidth}
        \centering
        \includegraphics[width=\textwidth]{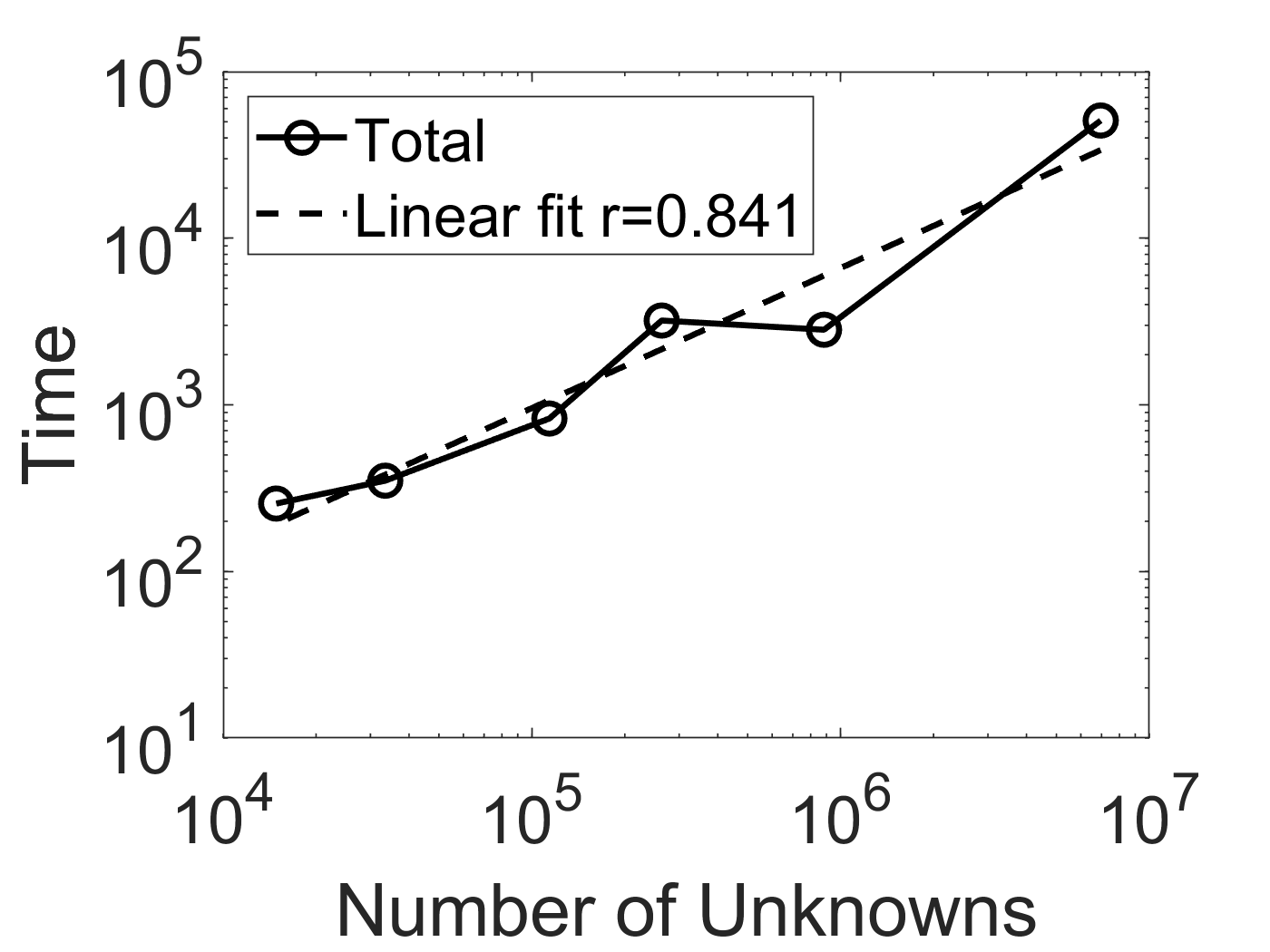}
        \caption{MSMS (100) GFM-ADI}
        \label{subfig:1a2e_spatial_time_gfmADI_msms100}
    \end{subfigure}
    \begin{subfigure}[b]{0.45\textwidth}
        \centering
        \includegraphics[width=\textwidth]{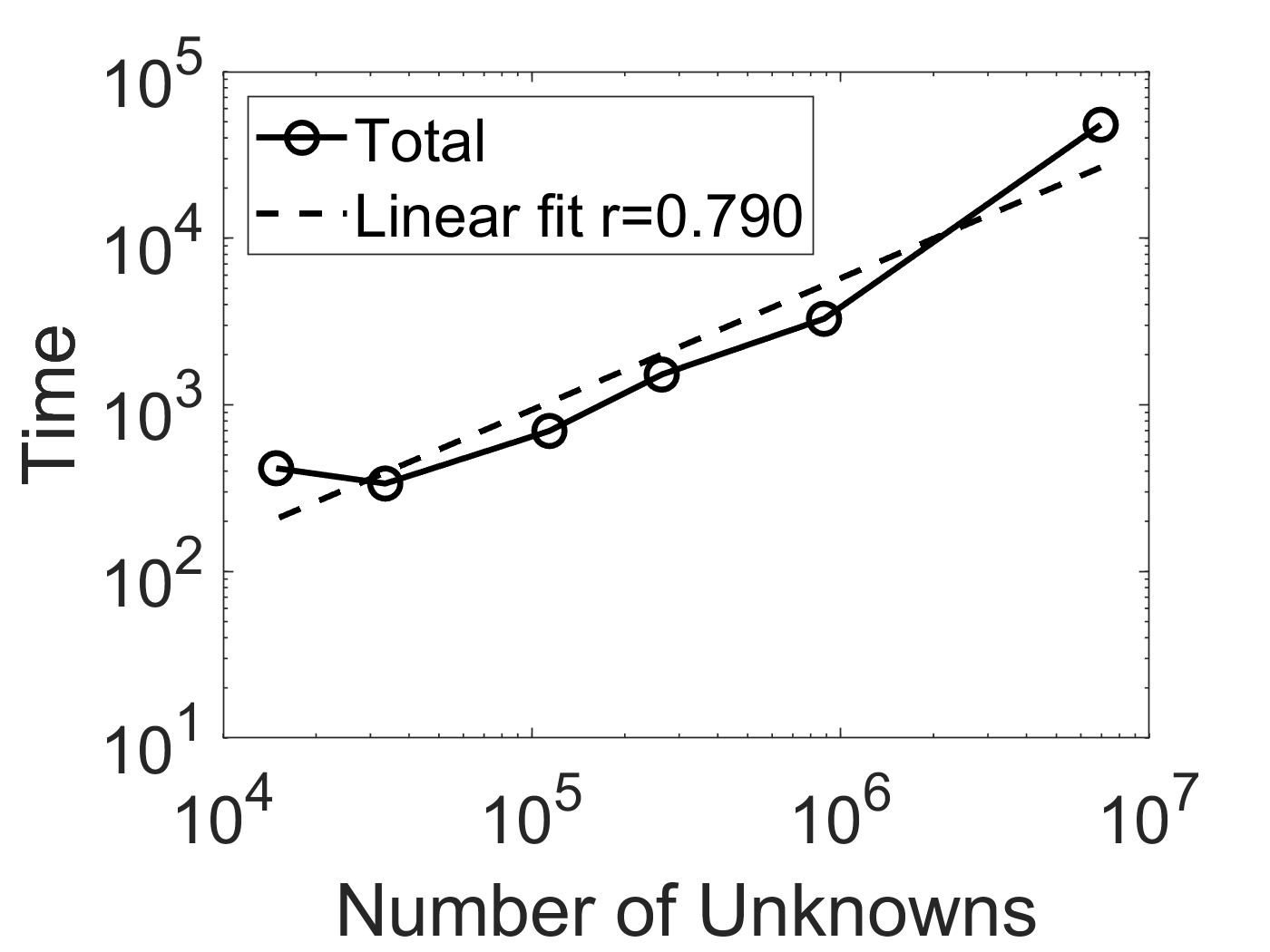}
        \caption{MSMS (100) GFM-LOD}
        \label{subfig:1a2e_spatial_time_gfmLOD_msms100}
    \end{subfigure}
    \caption{Total CPU time in seconds taken to compute the SES and solve the PBE for protein 1a2e. Also plotted are the least squares regression lines of best fit. }
    \label{fig:1a2e_spatial_times2}
\end{figure}

We next investigate the CPU time for protein 1a2e. The CPU time for solving the PBE and generating the SES is shown, respectively, in Fig. \ref{fig:1a2e_spatial_times} (a) and (b). In each subfigure, six curves are considered for two PBE solvers and three SES surfaces (ESES, MSMS with density 10, and MSMS with density 100). In Fig. \ref{fig:1a2e_spatial_times} (a), all six curves are similar. For part (b), it can be seen the MSMS(100) is more expensive than the MSMS(10) for all grid sizes. Here, the CPU time for the SES includes not only the SES generation, but also the Lagrangian-to-Eulerian conversion in case of the MSMS. For the MSMS generation, it takes the same CPU time for different grid size $h$, while the CPU for the conversion depends on the grid size $h$. In total, the CPU time for the MSMS barely increases or increases a little as $h$ becomes smaller. For the ESES, a different SES is generated at different $h$, so that the CPU time increases as total degree of freedom becomes larger. So the ESES is initially cheaper at a large $h$, and becomes as expensive as the MSMS(100) on the densest mesh. We also note that the CPU time for the GFM-ADI is slightly different from that of the GFM-LOD in Fig. \ref{fig:1a2e_spatial_times} (b), while they should be the same theoretically for generating molecular surface. The minor difference here is due to the randomness in calculating CPU time, which may be influenced by system processing time, memory reading time, etc.

The total CPU is shown in Fig. \ref{fig:1a2e_spatial_times} (c). Because the CPU time for the SES generation is much shorter than that for solving the PBE, the total CPU curves look like the part (a), and all six curves follow a similar trend. Based on the total CPU time, we next investigate the complexity of the pseudo-time solvers based on different SES definitions. For this purpose, we plot six cases separately in Fig. \ref{fig:1a2e_spatial_times2} against the total degree of freedom (DoF) $N$. The $\log$(CPU) as a function of $\log(N)$ is fitted with the least squares regression line. The slope $r$ of this line gives us the complexity order, i.e., the CPU is on the order of $O(N^r)$. It can be observed that in all cases, such rate is about one. This agrees with our finding in \cite{sheik_zhao_2020} that the complexity of the pseudo-time methods scale linearly with respect to the DoF $N$ when the number of time steps is fixed.  We note that the rates for the MSMS(100) are slightly less than the others. This does not mean the MSMS(100) is faster, but this is actually because the MSMS(100) takes more CPU time than the ESES and MSMS(10) for large $h$ values. 

\begin{table}[!htbp]
\centering
\begin{tabular}{|l|l|l|l|l|l|}\hline
PDB ID & \# Atoms & Method &ESES energy& MSMS energy& Relative difference \\\hline
\multirow{2}{*}{1bbl} & \multirow{2}{*}{576} & GFM-ADI & -984.01 & -986.53 & 2.55E-03\\
 & & GFM-LOD & -983.48 & -986.19 & 2.75E-03\\\hline
\multirow{2}{*}{1vii} & \multirow{2}{*}{596} & GFM-ADI & -898.29 & -901.78 & 3.87E-03 \\
 & & GFM-LOD & -897.99 & -901.37 & 3.75E-03 \\\hline 
  \multirow{2}{*}{1aho} & \multirow{2}{*}{962} & GFM-ADI & -887.46 & -893.86 & 7.16E-03 \\
 & & GFM-LOD & -884.42 & -890.76 & 7.11E-03 \\\hline
\multirow{2}{*}{1svr} & \multirow{2}{*}{1435} & GFM-ADI & -1704.58 & -1707.87 & 1.93E-03 \\
 & & GFM-LOD & -1703.00 & -1706.40 & 1.99E-03 \\\hline
\multirow{2}{*}{1a63} & \multirow{2}{*}{2065} & GFM-ADI & -2372.65 & -2370.80 & 7.79E-04 \\
 & & GFM-LOD & -2371.54 & -2369.29 & 9.50E-04 \\\hline
\multirow{2}{*}{1a7m} & \multirow{2}{*}{2809} & GFM-ADI & -2152.13 & -2155.05 & 1.36E-03 \\
 & & GFM-LOD & -2149.73 & -2152.50 & 1.29E-03 \\\hline
\end{tabular}
\caption{Comparison of the electrostatic free energies for six proteins 1bbl \cite{1BBL}, 1vii \cite{1VII}, 1aho \cite{1AHO}, 1svr \cite{1SVR}, 1a63 \cite{1A63}, and 1a7m \cite{1A7M} using the ESES and MSMS surfaces. MSMS used a density of $10$. Here, $\dt=0.005, h=0.5,$ and the stopping condition is either $T_{end}=10$ or $\Delta E_{sol}^n<10^{-4}$. The difference between reported energies relative to the MSMS energy is reported as the relative difference.}
\label{tab:5more}
\end{table}

We further validate the ESES by studying more proteins. The electrostatic free energies of six proteins are considered in Table \ref{tab:5more}. Because the instability issues associated with a high density, we only employed the density 10 for the MSMS. For the pseudo-time solvers, we take $h=0.5$ and $\dt=0.005$. The stopping condition is either $T_{end}=10$ or $\Delta E_{sol}^n<10^{-4}$. For each method, the energies calculated by the ESES and MSMS are reported and their relative difference is also shown. As demonstrated in the previous studies, the MSMS energy with density 10 is slightly inaccurate than that of the ESES. Thus, there is always some difference between energies of the ESES and MSMS. It can be observed in Table \ref{tab:5more}, such relative difference is not significant, ranging from 8E-04 to 7E-03. In general, the ESES energies are consistent with those of the MSMS. Because a better accuracy and stability, the present study validate the use of the ESES molecular surface in the pseudo-time GFM approach.

\begin{figure}[!htbp]
 \centering
 \begin{subfigure}[b]{0.45\textwidth}
 \centering
 \includegraphics[width=\textwidth]{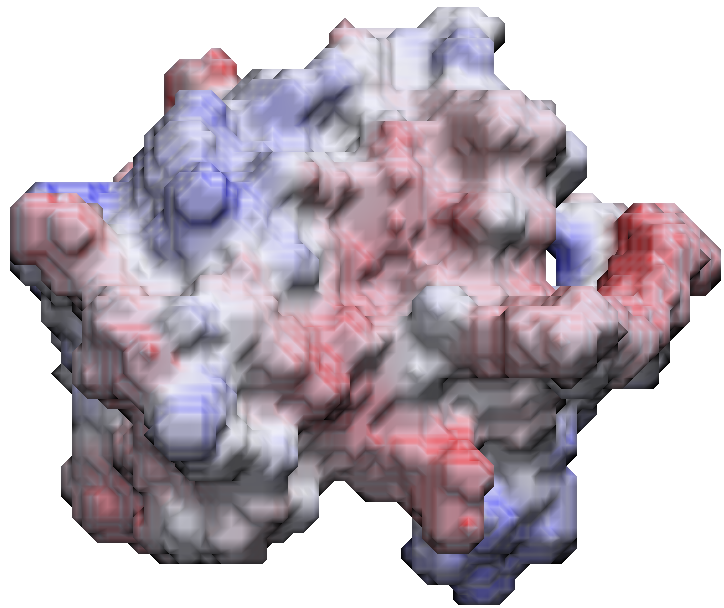}
 \caption{ESES and GFM-ADI}
 \label{subfig:1a2e_view_eses_gfmadi}
 \end{subfigure}
 \begin{subfigure}[b]{0.45\textwidth}
 \centering
 \includegraphics[width=\textwidth]{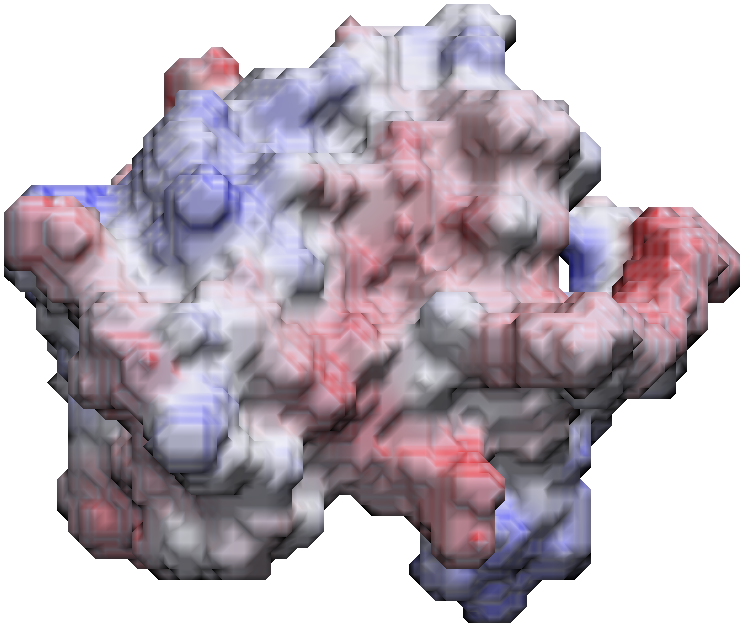}
 \caption{ESES and GFM-LOD}
 \label{subfig:1a2e_view_eses_gfmlod}
 \end{subfigure}
 \begin{subfigure}[b]{0.45\textwidth}
 \centering
 \includegraphics[width=\textwidth]{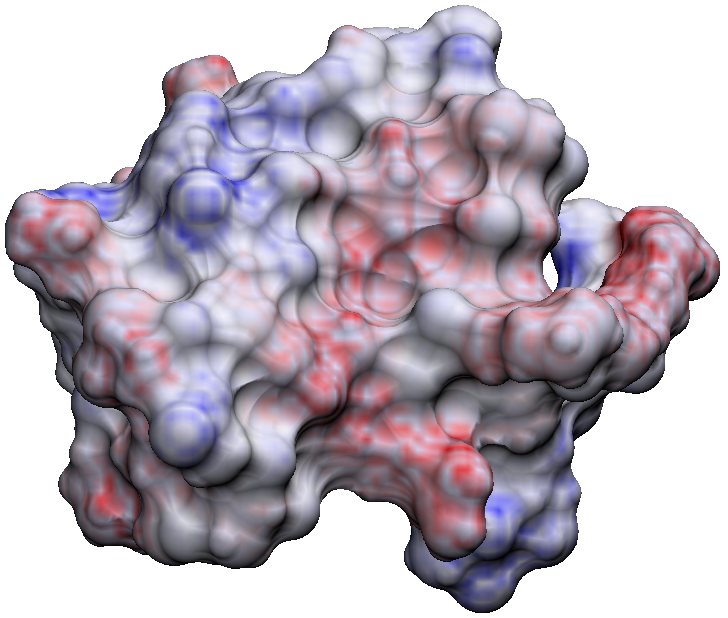}
 \caption{MSMS and GFM-ADI}
 \label{subfig:1a2e_view_msms_gfmadi}
 \end{subfigure}
 \begin{subfigure}[b]{0.45\textwidth}
 \centering
 \includegraphics[width=\textwidth]{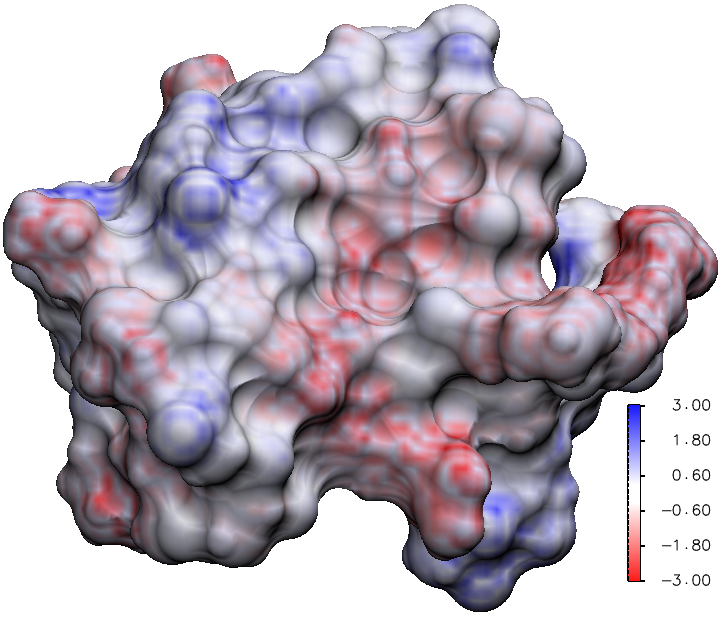}
 \caption{MSMS and GFM-LOD}
 \label{subfig:1a2e_view_msms_gfmlod}
 \end{subfigure}
 \caption{Potentials mapped onto molecular surface after solving the PBE for protein 1AHO with $h=0.5$, $\dt=0.005, T_{end}=10$. The MSMS surface used a density of $10$. Potentials were scaled linearly into the interval $[-3,3]$.}
 \label{fig:1a2e_view_potentials}
\end{figure}

We finally illustrate the visual difference between the ESES and MSMS, see Fig. \ref{fig:1a2e_view_potentials}. For each case, the potential calculated by the GFM-ADI or GFM-LOD has been mapped onto the molecular surface. It can be seen that two MSMS pictures are of high quality, because they are generated by using the VMD \cite{VMD} built-in functions. The ESES defines the same SES on Cartesian grids, but no suitable software is available for visualizing it. 
In the present study, an ESES hypersurface function is calculated as in Ref. \cite{ESES} and an isosurface is plotted by using the VMD. The ESES surfaces plotted in Fig. \ref{fig:1a2e_view_potentials} are clearly based on Cartesian volumetric data, and look slightly different from the MSMS surfaces. This visual resolution difference is because the visualization tool we used, VMD, supports better visualization of MSMS surfaces than the ESES surfaces. Note that the representations used in our calculations were based on Cartesian points for both the ESES and MSMS surfaces. This visualization resolution difference does not reflect any difference in surface quality for our numerical calculations.  Putting aside the quality difference due to visualization tools, the four surfaces in Fig. \ref{fig:1a2e_view_potentials} are actually very similar, in terms of both the shape and the mapped potentials.

\section{Adaptive pseudo-time GFM approach}
In this section, we explore hypothesis and techniques for adaptive time integration in the pseudo-time GFM simulations. In all tests, the nonlinear PBE is solved with $h=0.5$, $\epsilon^+=80$, $\epsilon^-=1$, and $\kappa =0.15$ M, and the ESES molecular surface is employed. Unless specified otherwise, the initial condition is taken as  the solution of the linearized PBE.

\begin{table}[!htbp]
\centering
\begin{tabular}{|l|l|l|l|l|l|l|l|}\hline
Method    & Test \# & $E_{sol}$    & CPU (s)    & Initial $\dt$ & final $\dt$ & Switch time & Error      \\\hline
\multirow{6}{*}{GFM-ADI} 
& 1 & -1140.343 & 3819.01 & 0.001     & 0.001   & N/A        & N/A        \\
& 2 & -1138.156 & 292.14  & 0.01      & 0.01   & N/A        & 0.00191785  \\
 &3  & -1139.766 & 522.61  & 0.01      & 0.001   & 9          & 0.00050631  \\
& 4 & -1138.185 & 519.98  & 0.001     & 0.01    & 1          & 0.00189262  \\
& 5 & -1140.205 & 1421.55 & 0.01      & 0.001   & 5          & 0.00012112  \\
& 6& -1138.295 & 1424.10 & 0.001     & 0.01    & 5          & 0.00179579 \\\hline
\multirow{6}{*}{GFM-LOD} 
& 1& -1140.581 & 5018.96 & 0.001     & 0.001   & N/A        & N/A        \\
& 2 & -1139.139 & 336.98  & 0.01      & 0.01   & N/A        & 0.00126422 \\
&3  & -1140.333 & 613.13  & 0.01      & 0.001   & 9          & 0.00021708 \\
&4 & -1139.147 & 620.03  & 0.001     & 0.01    & 1          & 0.00125672 \\
&5  & -1140.537 & 1673.64 & 0.01      & 0.001   & 5          & 0.00003841 \\
 &6 & -1139.183 & 1653.28 & 0.001     & 0.01    & 5          & 0.00122576 \\\hline
\end{tabular}
\caption{The electrostatic free energies calculated by using different combinations of $\dt$ values
and the corresponding CPU time. The protein 2pde is used with $h=0.5$ and $T_{end}=10$. }
\label{tab:hypothesis}
\end{table}

\subsection{Hypothesis test}
The adaptive time stepping has been widely studied in the literature \cite{Bank1980,Kelley1998,Pollock2015} for an efficient convergence to the steady state. Because the solution approaches a limit as the time $t$ becomes larger, the temporal variation of the solution becomes smaller. Thus, a common hypothesis in the field is that as $t$ increases, a larger $\dt$ is sufficient to capture temporal variation, which saves total time steps for accelerating the convergence. In the present study, we first test this hypothesis by considering a protein 2pde \cite{2PDE} under different combinations of $\dt$ values. For simplicity, computations were stopped when $T_{end}=10$. For both GFM-ADI and GFM-LOD, 6 tests were carried out. In Test 1 and Test 2, a constant $\dt=0.001$ and $\dt=0.01$ were employed, respectively. Test 3 used $\dt=0.01$ until $t=9$, then $\dt=0.001$. Test 4 used $\dt=0.001$ until $t=1$, then $\dt=0.01$. Test 5 used $\dt=0.01$ until $t=5$, then $\dt=0.001$. Test 6 used $\dt=0.001$ until $t=5$, then $\dt=0.01$. Note that each value of $\dt$ was used for the same number of steps in tests 3 and 4. The same is true for tests 5 and 6. The resulting energies and the corresponding CPU time are listed in Table \ref{tab:hypothesis}. For both the GFM-ADI and GFM-LOD schemes, the energy of Test 1 is chosen as the reference value, and the relative errors for other tests are also reported. 

By examining the results in Table \ref{tab:hypothesis}, the following conclusion can be drawn, which is true for both the GFM-ADI and GFM-LOD schemes. By comparing Test 3 with Test 4 and comparing Test 5 with Test 6, the test with the smaller final value of $\dt$ is always more accurate, while both tests consume similar execution times because of the same number of time steps. In section \ref{large-study}, this finding has been further verified in many other proteins for the present pseudo-time GFM approach. This means that the common hypothesis is invalid for the pseudo-time PBE. 
The reason for this is perhaps because of the PBE nonlinearity and its treatment. A time splitting is the key of the pseudo-time method for treating PBE nonlinearity analytically, which inevitably introduces a splitting error on the order of $O(\dt)$. Such a splitting error can be controlled only if a small enough $\dt$ is used, even when the steady-state is approaching. 

Moreover, if one compares among Test 2, 4, and 6, i.e., three tests with the same ending $\dt$, one can see their accuracies are almost the same. In other words, the use of a smaller $\dt$ in early stage does not improve the final accuracy. It is the final $\dt$ value that determines the final accuracy. Thus, the use a smaller $\dt$ in early stage should be avoided, because it simply wastes the CPU time. Therefore, in order to save computation time without a significant loss of accuracy, the best strategy is to start with a large $\dt$ and reduce $\dt$ as $t$ increases. This is the basic idea for our further exploration. 

\begin{table}[!htbp]
\centering
\begin{tabular}{|l|l|l|r|l|l|l|}\hline
Method           & $\dt_{max}$ & $\dt_{min}$ & CPU (s) & Energy    & Relative CPU & Relative error \\\hline
\multirow{4}{*}{Constant $\dt$} 
       & 0.001     & 0.001     & 17535   & -1157.652 & N/A     & N/A      \\
        & 0.010     & 0.010     & 1944    & -1155.439 & 11.09\%     & 1.912E-03       \\
         & 0.100     & 0.100     & 215     & -1143.129 & 1.223\%     & 1.254E-02      \\
         & 1.000     & 1.000     & 40      & -1112.360 &  0.229\%     & 3.912E-02      \\\hline
MANUAL1          & 1.000     & 0.001     & 9967    & -1156.932 & 56.84\%     & 6.218E-04      \\\hline
MANUAL2          & 1.000     & 0.001     & 5505    & -1156.945 & 31.40\%     & 6.109E-04      \\\hline
PID1      & 1.000     & 0.001     & 3850    & -1157.038 & 21.95\%     & 5.303E-04      \\\hline
PID2      & 1.000     & 0.001     & 9598    & -1157.651 & 54.74\%     & 3.790E-07      \\\hline
FastPID          & 1.000     & 0.001     & 95      & -1145.151 & 0.544\%    & 1.080E-02      \\\hline
NonincreasingPID & 1.000     & 0.001     & 4840    & -1156.900 & 27.60\%     & 6.494E-04      \\\hline
\end{tabular}
\caption{The electrostatic free energies and CPU time of the protein 2pde for different constant $\dt$ values and various adaptive $\dt$ schemes. The relative CPU and error are calculated by using the result of $\dt=0.001$ as the reference solution.}
\label{tab:2pde}
\end{table}

\subsection{One protein study}
For the rest of this paper, we will consider the GFM-LOD method only. The developed adaptive schemes can be applied to the GFM-ADI too. But the use of the GFM-LOD gives us more freedom to search for the best adaptive strategy, because any $\dt$ value can be employed, thanks to the unconditional stability of the GFM-LOD \cite{sheik_zhao_2020}. 

We continue to study the protein 2pde in this subsection by designing various time stepping schemes. To benchmark different schemes, a reference solution is generated with $\dt = 0.001$ and stopping condition being either $T_{end}=50$ or $\Delta E^n_{sol}< 10^{-4}$. See Table \ref{tab:2pde}. We note that such a computation is very costly, and is not recommended for practical usage. The results by using $\dt=0.01$, $0.1$ and $1$ are also listed in Table \ref{tab:2pde}. It can be seen that as $\dt$ increases, the relative error becomes larger, while the CPU time becomes smaller. These constant $\dt$ results also provide benchmarks for examining adaptive $\dt$ schemes. 

We consider a series of methods with $\dt_{max}=1$ and $\dt_{min}=0.001$, and allow $\dt$ to be changed only in the range of $[\dt_{min},\dt_{max}]$, i.e., $\dt_{min} \le \dt \le \dt_{max}$. We first explore two ``manual'' methods, in which the refinement of $\dt$ is automatically applied whenever a criterion in terms of an error measurement $e$ is satisfied. In particular, two error norms are considered to measure the temporal variation, i.e., the change in two steps for the solution $e=|| U^n - U^{n-1}||_2$ and that for the free energy $e=\Delta E^n_{sol}$. The corresponding methods are called Manual 1 and Manual 2, respectively. In both methods, we take $\dt = \dt_{max}$ and $\delta=1$ initially. In each time step, we keep checking if $e < \delta$. When $e < \delta$, we halve both $\dt$ and $\delta$, i.e., $\dt=\dt/2$ and $\delta=\delta/2$. After the condition $e < \delta$ being met for multiple times, the minimal $\dt$ is reached and no further refinement is conducted whenever $\dt=\dt_{min}$. The computation is continued until either $T_{end}=50$ or $\Delta E^n_{sol} < 10^{-4}$. It can be seen from Table \ref{tab:2pde} that both manual methods produce very good accuracy, while save CPU time. Moreover, the Manual 2 scheme is better, because it demands less CPU time. 

The time increment $\dt$ produced by two manual methods is essentially a piecewise constant, i.e., it is a constant for a certain time period and is halved for the next period. This is not an adaptive $\dt$ process. The adaptive time integration in numerical analysis means that one keeps tracking the temporal variation or error and based on that to select $\dt$ simultaneously. For this purpose, we adopt the PID control technique developed in \cite{pid} to track the temporal change and correspondingly define a scaling factor
\begin{equation}\label{F-factor}
	F=\left(\frac{e_{n-1}}{e_n}\right)^{k_P}\left(\frac{\epsilon_p}{e_n}\right)^{k_I}\left(\frac{e_{n-1}^2}{e_n e_{n-2}}\right)^{k_D},
\end{equation}
where the same parameters as in \cite{pid} are used, i.e.,  $\epsilon_p=0.0025, k_P=0.075, k_I=0.175,$ and $k_D=0.01$. In each time step, the time increment is updated as $\dt=\frac{\dt}{F}$. Here $e_n$ is an error norm. In the present study, we test two PID schemes with $e_n$ being chosen as
\begin{align}
    e^{(u)}_n &= \frac{||u^n-u^{n-1}||_2}{||u^n||_2}, \\
    e^{(E)}_n &= \left|\frac{\Delta E_{sol}^n}{E_{sol}^n}\right|, 
\end{align} 
for the PID1 and PID2, respectively. In our computations, the factor $F$ is bounded between $0.2$ and $5.0$, and $\dt_{min} \le \dt \le \dt_{max}$ is always ensured. The time integration will be stopped when either $T_{end}=50$ or $\Delta E_{sol}^n <10^{-4}$. The key features of the PID1 and PID2 are summarized in Table \ref{tab:PID}. 

\begin{table}[!htbp]
\centering
\begin{tabular}{|l|l|l|l|}\hline
Name             & $e_n $      & Stopping condition \\\hline
PID1      & $e_n^{(u)}$  & $T_{end}=50$ or $\Delta E_{sol}^n <10^{-4}$                                                    \\\hline
PID2      & $e_n^{(E)}$  & $T_{end}=50$ or $\Delta E_{sol}^n < 10^{-4}$                                                    \\\hline
FastPID          & $e_n^{(u)}$  & $T_{end}=50$ or $\Delta E_{sol}^n < 10^{-4}$ or 100 steps after   $\dt_{min}$ is first reached. \\\hline
NonincreasingPID & $e_n^{(u)}$  & $T_{end}=50$ or ($\Delta E_{sol}^n < 10^{-4}$ and $\dt_{min}$ is reached).\\\hline
\end{tabular}
\caption{PID methods with different error norms and stopping conditions. The NonincreasingPID method has an additional constraint to prevent $\dt$ from increasing at any time step.}
\label{tab:PID}
\end{table}

The numerical results of the PID1 and PID2 are also listed in Table \ref{tab:2pde}. The PID2 has an extremely small error, while the accuracy of the PID1 is better than two manual methods. On the other hand, the PID2 is quite expensive, while the PID1 is faster than the PID2 and two manual methods. We note that the adaptive time process essentially trades the accuracy for efficiency. It is not necessary to achieve a high precision like the PID2. Instead, we aim to maintain the relative error less than 1.0E-3, while reduce the CPU time as much as possible. Based on these considerations, the PID1 is obviously our pick, among two PID and two manual methods. 

Two more PID methods, i.e., FastPID and NonincreasingPID, are also tested in Table \ref{tab:2pde}. For the protein 2pde, these two PID methods perform worse than the PID1 -- they were designed after more protein studies have been conducted. The two new PID methods are designed based on the PID1, i.e., using the same error norm $e_n = e_n^{(u)}$, but have different stopping condition, see Table \ref{tab:PID}. For the FastPID, the computation could stop earlier after $\dt$ has been reduced to $\dt=\dt_{min}$ and then 100 steps further computations are conducted. As its name stands, the FastPID is very efficient, but its error is about 1\% for the protein 2pde. The original PID method \cite{pid} allows $\dt$ to be changed according to the temporal variation, which could result in a larger or smaller new $\dt$. As explained above, our aim is to reduce $\dt$ monotonically in the pseudo-time GFM approach. A NonincreasingPID is thus designed, in which the scaling factor $F$ calculated from (\ref{F-factor}) will be reset to be one when $F < 1$. This makes sure that $\dt$ decreases monotonically in the NonincreasingPID. 

\subsection{Large scale study}\label{large-study}
In this subsection, we test the PID methods by considering a test set of 74 proteins studied in \cite{proteins}, which gives a representative sample of proteins in the protein databank \cite{RCSB_PDB}. We refer to \cite{proteins} for more details on how this set is selected and processed based on specific physical features. A reference solution is generated for each protein in this set with a constant $\dt=0.01$ and stopping condition being either $T_{end}=50$ or $\Delta E^n_{sol} < 10^{-4}$. Note that the constant $\dt$ is larger in this large scale study in order to save simulation time for a total of 74 proteins. Moreover, a minimal stopping time is imposed so that the computation will not stop before $t=5$ even if $\Delta E^n_{sol} < 10^{-4}$ is satisfied. Details about this early stopping issue will be discussed in next subsection.

\begin{table}[!htbp]
\centering
\begin{tabular}{|l|l|l|l|}\hline
$\dt_{max}$         & Maximum relative error & Mean relative error & Mean relative time \\\hline
1.0 & 0.1037399   & 0.0485009               & 0.135525608           \\\hline
0.5 & 0.0666789   & 0.0172739               & 0.163347295           \\\hline
0.1 & 0.0151291   & 0.0033452               & 0.161796101           \\\hline

\end{tabular}
\caption{Effect of $\dt_{max}$ on the accuracy and efficiency of the FastPID for 74 proteins.  }
\label{tab:fastpid}
\end{table}

We first test the performance of the FastPID by using $\dt_{max}=1$ and $\dt_{min}=0.001$. The FastPID is almost the same as the PID1, except an early stop is allowed after 100 steps of integration are conducted at $\dt=\dt_{min}$. The FastPID prevents time-stepping at $\dt=\dt_{min}$ for a long time in the PID1, which will be very expensive. 

\begin{figure}[!htbp]
    \centering
    \includegraphics[width=\textwidth]{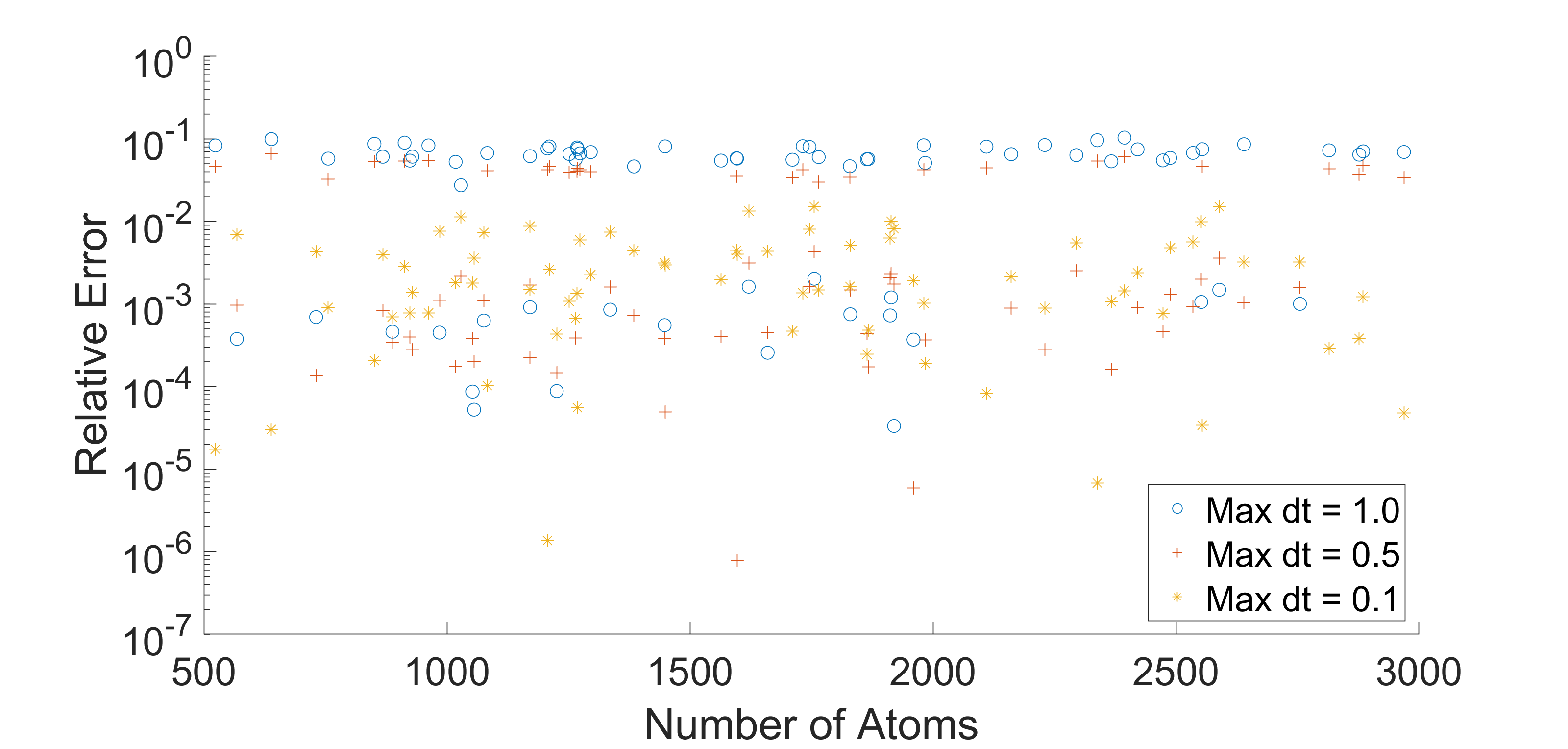}
    \caption{Relative errors in electrostatic free energy for 74 proteins computed with the FastPID method using $\dt_{max}=1.0, 0.5,$ and $0.1$.}
    \label{fig:error-fastpid-maxdts}
\end{figure}

The FastPID results are reported in Table \ref{tab:fastpid} for different $\dt_{max}$. By using $\dt_{max}=1$, the maximum relative error is about 10\%, which is intolerable for real applications. By checking different parameters of the FastPID, it is found in our studies that the max error could be reduced by using a smaller $\dt_{max}$. It can be seen in Table \ref{tab:fastpid} that the max relative error becomes 7\% and 2\%, respectively, for $\dt_{max}=0.5$ and $\dt_{max}=0.1$. In order to help us to comprehend these results, we plot the relative errors of the FastPID method for 74 proteins in Fig. \ref{fig:error-fastpid-maxdts}. We see that there is a strip of errors between $10^{-2}$ and $10^{-1}$ with $\dt_{max}=1.0$ and $\dt_{max}=0.5$, but not $\dt_{max}=0.1$. In particular, there are many proteins whose relative errors are large, not just a few of them. This study suggests that the FastPID scheme fails to be a good adaptive time integration technique, because the final accuracy depends on the initial $\dt=\dt_{max}$. If the adaptive process was effective, it is the final $\dt=\dt_{min}$ that determines the final accuracy. Then the accuracy should not critically depend on $\dt_{max}$.

\begin{figure}[!htbp]
    \centering
 \begin{subfigure}[b]{0.75\textwidth}
 \centering
 \includegraphics[width=\textwidth]{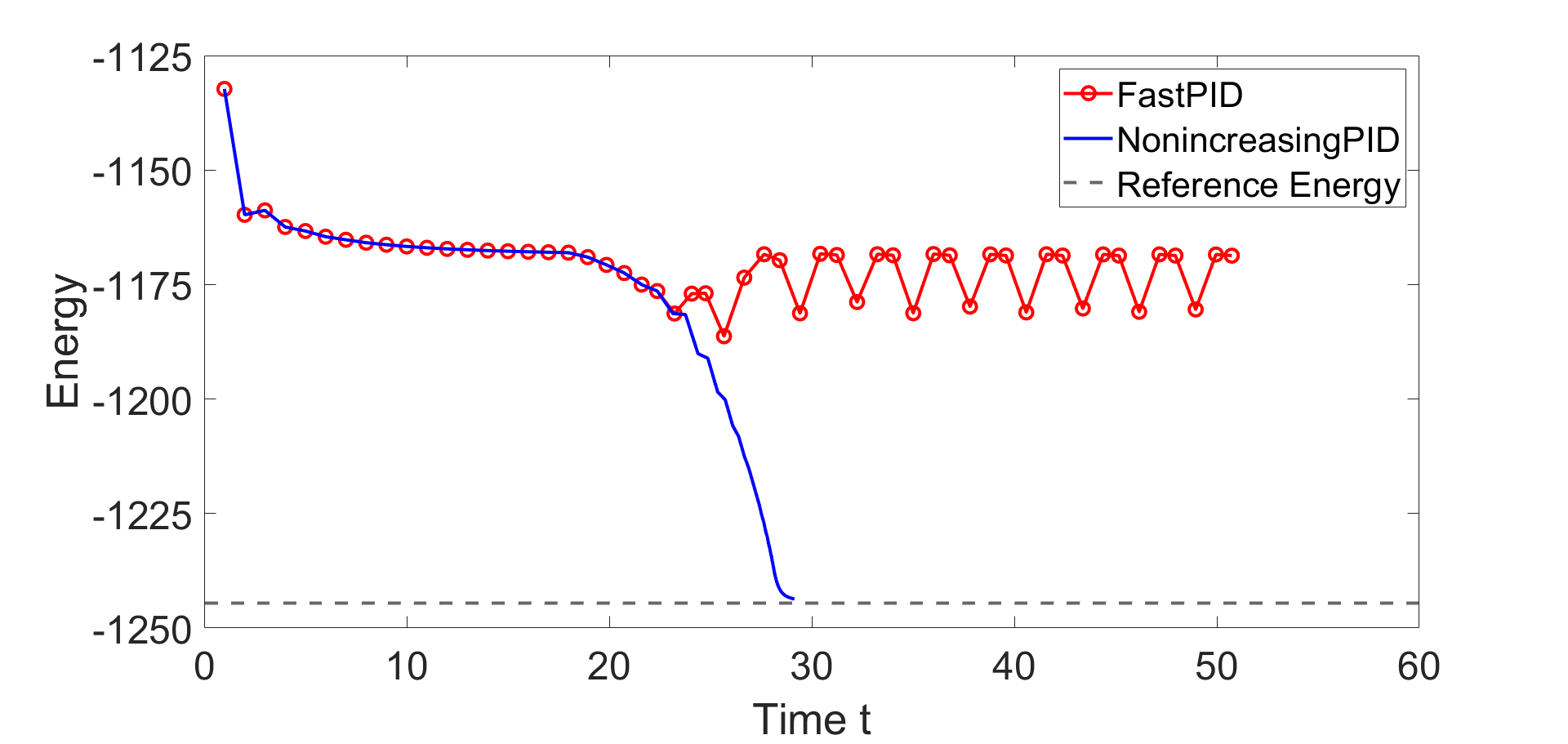}
 \caption{Electrostatic free energy}
\label{fig:oscillation-energy}
 \end{subfigure}
  \begin{subfigure}[b]{0.75\textwidth}
 \centering
 \includegraphics[width=\textwidth]{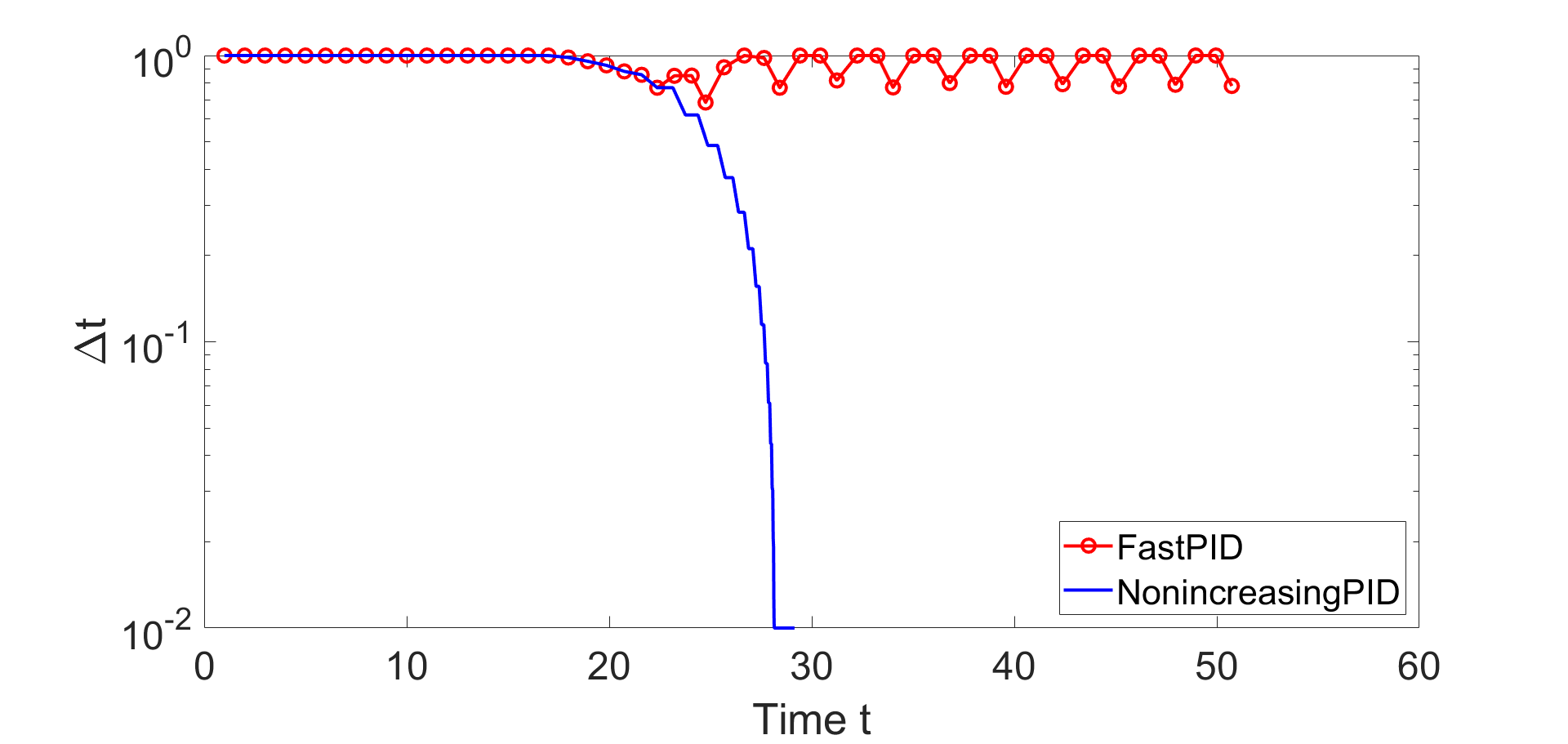}
 \caption{$\dt$}
\label{fig:oscillation-dt}
 \end{subfigure}
   \caption{The electrostatic free energy and $\dt$ value over time for the protein 1zuu. For the FastPID, we take $\dt_{max}=1.0$. The reference energy is plotted as a horizontal line in part (a). }
    \label{fig:1zuu}
\end{figure}
In searching for the cause for the failure of the FastPID scheme, we found that for many proteins, the FastPID never reaches $\dt=\dt_{min}$ in the PID updates. One such example is shown in Fig. \ref{fig:1zuu}. By plotting the FastPID free energy $E_{sol}$ and $\dt$ against the time $t$, we see that both of them oscillate for the second half of the time integration period until the computation is stopped by $T_{end}=50$. In particular, after $\dt$ becomes smaller than $\dt_{max}$ around $t=20$, the PID calculation generates a factor $F <1$, so that $\dt$ becomes larger and bounces back to $\dt=\dt_{max}$. There is a competing process so that $\dt$ oscillates just below $\dt_{max}=1$. Consequently, the energy is trapped in the oscillation, and never converges to the steady state. On the other hand, we note that in the FastPID cases where an oscillation occurs, $\dt$ never reaches $\dt_{min}$. Then PID1 and FastPID will behave identically in these cases. This is why we did not report results of the PID1 -- they are also affected by the $\dt$ oscillation issue. 

By reducing the value of $\dt_{max}$, it is possible that the simulation will oscillate around the new $\dt_{max}$. It is also possible that the adaptive process could overcome the barrier and converge to the right place. This is essentially why the error of the FastPID depends on $\dt_{max}$ in Fig. \ref{fig:error-fastpid-maxdts}. However, even with $\dt_{max}=0.1$, the max relative error is still $0.015$ in Table \ref{tab:fastpid}. This means that one cannot count on changing parameter values to suppress the $\dt$ oscillation issue. A systematic change has to be conducted. 

The NonincreasingPID method is designed to overcome the $\dt$ oscillation issue of the PID1 and FastPID schemes. For the PID scaling factor $F$ calculated from (\ref{F-factor}), it will be reset to be one when $F < 1$. This makes sure that $\dt$ decreases monotonically in the NonincreasingPID. It can be seen in Fig. \ref{fig:1zuu} that the NonincreasingPID energy is the same as that of the FastPID initially, but later converges correctly to the reference solution. Changing monotonically, $\dt$ reaches $\dt_{min}$ in the NonincreasingPID. Shortly after that, the computation is stopped due to $E^n_{sol} < 10^{-4}$ before $T_{end}=50$. 

\begin{table}[!htbp]
\centering
\begin{tabular}{|l|l|l|l|l|}\hline
$\dt_{max}$                  & $\dt_{min}$                  & $TOL$ & Mean relative error & Mean relative time \\\hline
\multirow{2}{*}{1.000} & \multirow{2}{*}{0.100} & 0.10      & 0.02502945          & 0.1289             \\
                       &                        & 0.01      & 0.00886432          & 0.1422             \\\hline
\multirow{2}{*}{1.000} & \multirow{2}{*}{0.010} & 0.10      & 0.00515138          & 0.1307             \\
                       &                        & 0.01      & 0.00060521          & 0.1618             \\\hline
\multirow{2}{*}{1.000} & \multirow{2}{*}{0.001} & 0.10      & 0.02240487          & 0.1303             \\
                       &                        & 0.01      & 0.02082860          & 0.2009             \\\hline
\multirow{2}{*}{0.100} & \multirow{2}{*}{0.010} & 0.10      & 0.00424923          & 0.1297             \\
                       &                        & 0.01      & 0.00208119          & 0.1855             \\\hline
\multirow{2}{*}{0.100} & \multirow{2}{*}{0.001} & 0.10      & 0.00592353          & 0.1194             \\
                       &                        & 0.01      & 0.00404972          & 0.2194             \\\hline
\multirow{2}{*}{0.010} & \multirow{2}{*}{0.001} & 0.10      & 0.00306056          & 0.3159             \\
                       &                        & 0.01      & 0.00182219          & 0.3407            \\\hline
\end{tabular}
\caption{Mean relative error in electrostatic free energy 
and mean relative CPU time of the NonincreasingPID using different $\dt_{max}, \dt_{min}$, and TOL. }
\label{tab:nonincreasingPID_params}
\end{table}

To optimize the performance of the NonincreasingPID scheme, we explore different parameter values. The stopping condition of the NonincreasingPID is still $T_{end}=50$ or ($\Delta E_{sol}^n <$ TOL and $\dt_{min}$ is reached), while now we take TOL as a parameter. We will also change $\dt_{max}$ and $\dt_{min}$ for $\dt_{min} \le \dt \le \dt_{max}$. The final $\dt=\dt_{min}$ controls the accuracy. The span between $\dt_{max}$ and $\dt_{min}$, as well as TOL, adjusts the length of time integration. Various combinations of these three parameters are studied and their average results are reported in Table \ref{tab:nonincreasingPID_params}. 
It can be observed that setting $\dt_{max}=1.0, \dt_{min}=0.01,$ and $TOL=0.01$ yields an extremely accurate and efficient solution with a mean relative error of $0.06\%$ and average relative computation time of $16.18\%$. Table \ref{tab:nonincreasingPID_params} shows that substantially changing any of the parameters produces a significant amount of error. Moreover, those that are more efficient offer only small time savings in exchange for substantial accuracy losses. The NonincreasingPID scheme with this set of parameters is our recommendation for future developments of the pseudo-time GFM approach.

\subsection{Initial condition}
We finally consider an important issue related to the pseudo-time integration, i.e., the initial condition. In the pseudo-transient approach to the nonlinear PBE, the initial condition is usually chosen either as a zero solution or the numerical solution of the linearized PBE (LPB) \cite{Geng13,Zhao14}. Theoretically, if the steady state solution exists with an appropriate initial condition, it satisfies the original boundary value problem (BVP) of the PBE.  Since the solution to the original BVP is known to be unique, the steady state solution is also unique. In the pseudo-time GFM approach, it have been observed \cite{sheik_zhao_2020} that the electrostatic free energies based on both zero and LPB initial solutions always converge to the same value after a long time integration. It is interesting to explore this for the NonincreasingPID scheme and study the related issues.  

\begin{figure}
    \centering
    \begin{subfigure}[b]{.45\textwidth}
    \includegraphics[width=\textwidth]{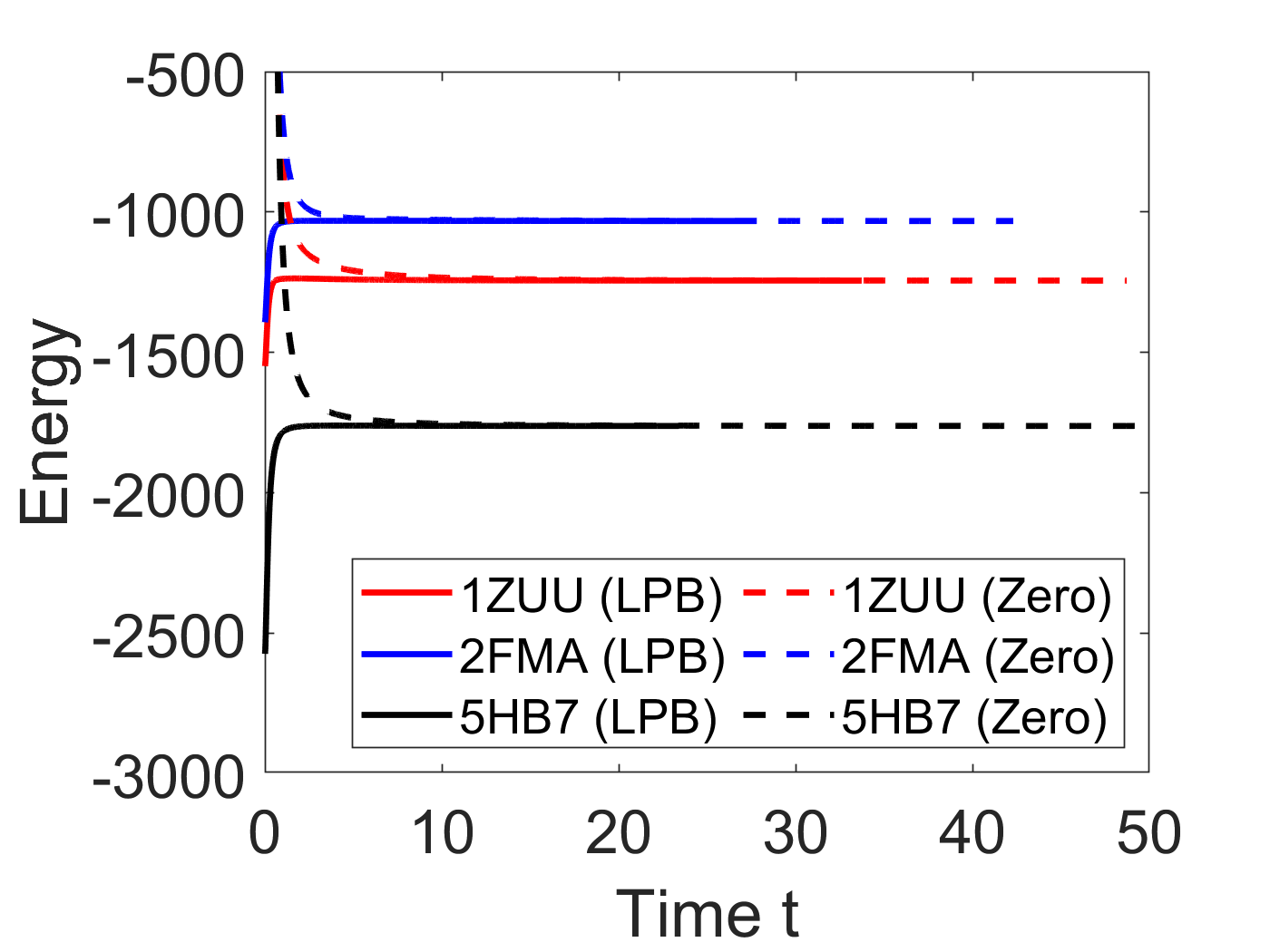}
    \caption{Constant $\dt$}
    \label{subfig:initialcond-energy-const}
    \end{subfigure}
    \begin{subfigure}[b]{.45\textwidth}
    \includegraphics[width=\textwidth]{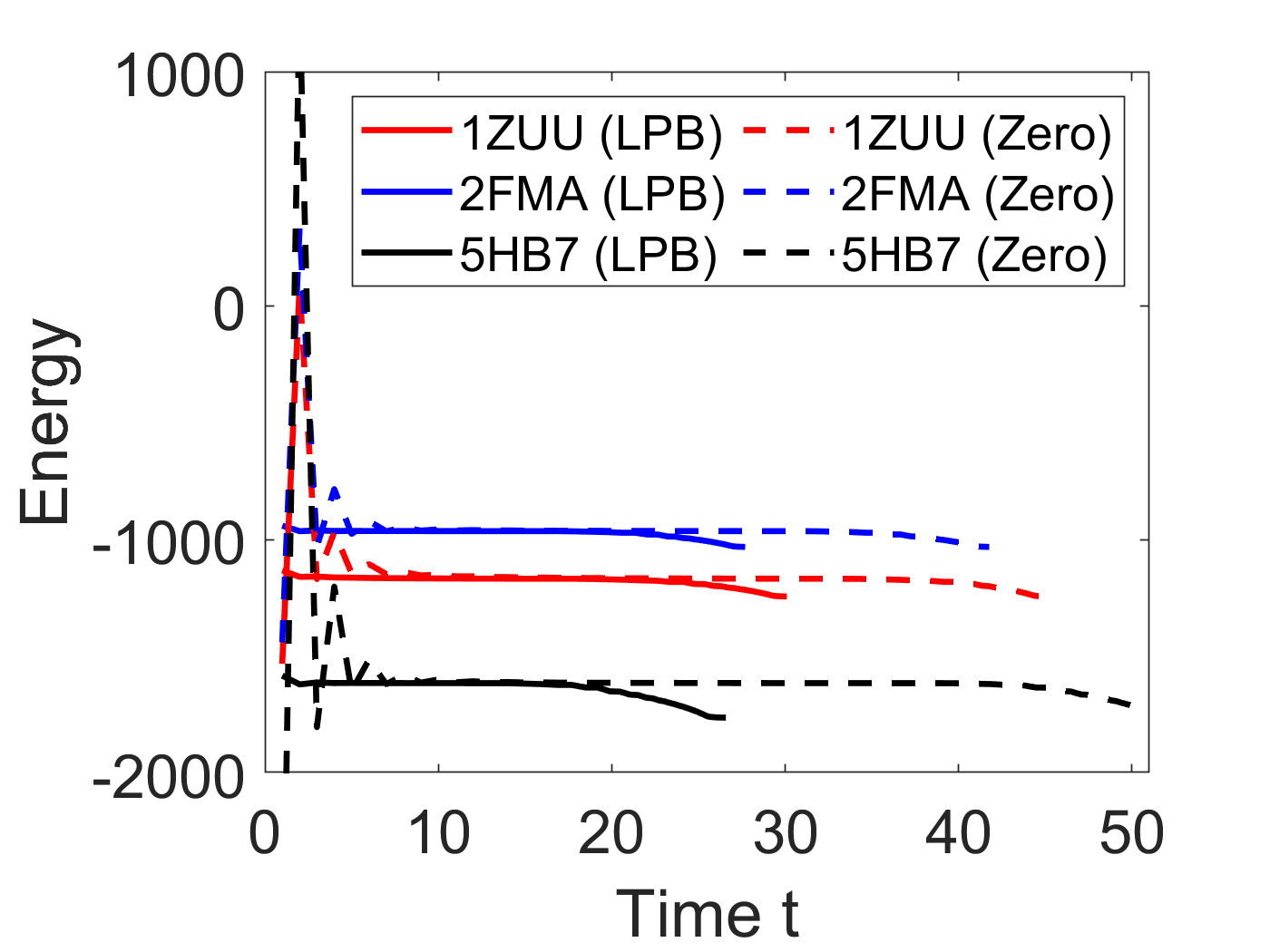}
    \caption{NonincreasingPID}
    \label{subfig:initialcond-energy-nonincreasing}
    \end{subfigure}
    \caption{Electrostatic free energy vs time for proteins 1ZUU, 2FMA, and 5HB7 solved with the zero and LPB initial conditions. The method for calculating  $\dt$ is constant in (a) and NonincreasingPID in (b).}
    \label{fig:initialcond-energy}
\end{figure}

\begin{table}[!htbp]
\centering
\begin{tabular}{|l|l|l|l|l|}\hline
protein               & method                            & initial condition & relative error & relative time \\\hline
\multirow{4}{*}{1ZUU} & \multirow{2}{*}{constant $\dt=0.01$}         & LPB               & N/A            & N/A           \\
                      &                                   & Zero              & 7.97844E-11    & 1.343608410    \\\cline{2-5}
                      & \multirow{2}{*}{NonincreasingPID} & LPB               & 0.000689049    & 0.146170049   \\
                      &                                   & Zero              & 0.000685721    & 0.092432609   \\\hline
\multirow{4}{*}{2FMA} & \multirow{2}{*}{constant $\dt=0.01$}         & LPB               & N/A            & N/A           \\
                      &                                   & Zero              & 1.88456E-10    & 1.229137136   \\\cline{2-5}
                      & \multirow{2}{*}{NonincreasingPID} & LPB               & 0.000445416    & 0.088169886   \\
                      &                                   & Zero              & 0.000442035    & 0.157473301   \\\hline
\multirow{4}{*}{5HB7} & \multirow{2}{*}{constant $\dt=0.01$}         & LPB               & N/A            & N/A           \\
                      &                                   & Zero              & 2.15751E-10    & 0.832033498   \\\cline{2-5}
                      & \multirow{2}{*}{NonincreasingPID} & LPB               & 0.000378474    & 0.137279424   \\
                      &                                   & Zero              & 0.000404068    & 0.142998219   \\\hline
\end{tabular}
\caption{Relative errors in electrostatic free energies and relative computation times for proteins 1ZUU, 2FMA, and 5HB7. For each protein, the reference solution is generated by the constant $\dt$ method with the LPB initial condition. Note that the NonincreasingPID calculation with the zero initial condition for protein 5HB7 used a larger $T_{end}=75$ because the value of $\dt$ did not reach $\dt_{min}$ before $t=50.$ }
\label{tab:initialcond}
\end{table}

By considering three proteins, 1ZUU, 2FMA, and 5HB7, the electrostatic free energies calculated based on the zero and LPB initial conditions are plotted against the time $t$ in Fig. \ref{fig:initialcond-energy} for both constant $\dt$ and NonincreasingPID schemes. For the constant $\dt$ method, the uniqueness is obvious. Two energy curves of a protein always approaches to the same limit, with one from above and one from below. On the other hand, the figure for the NonincreasingPID scheme looks different. Two curves attains the same intermediate energy for certain time period, but each curve makes a drop later. This is because in the NonincreasingPID scheme, a larger $\dt > \dt_{min}$ is operated in the middle course of the process. Once $\dt = \dt_{min}$, the energy quickly converges to the true answer, which is below the intermediate value. So the drops in the Fig. \ref{fig:initialcond-energy} (b) represent the final convergence stage. 

The converged values of the zero and LPB initial conditions are illustrated in Table \ref{tab:initialcond}. For the constant $\dt$ method, the relative energy difference between the zero and LPB initial conditions is as small as $10^{-10}$. For the NonincreasingPID scheme, both energies agree with that of the constant $\dt$ very well, and the relative energy difference between the zero and LPB initial conditions is about $10^{-5}$. This quantitatively verified that even though they appear different in the Fig. \ref{fig:initialcond-energy} (b), the zero and LPB initial conditions actually produce the same steady state solution for the NonincreasingPID scheme.

\begin{figure}
    \centering
    \begin{subfigure}[b]{.45\textwidth}
    \includegraphics[width=\textwidth]{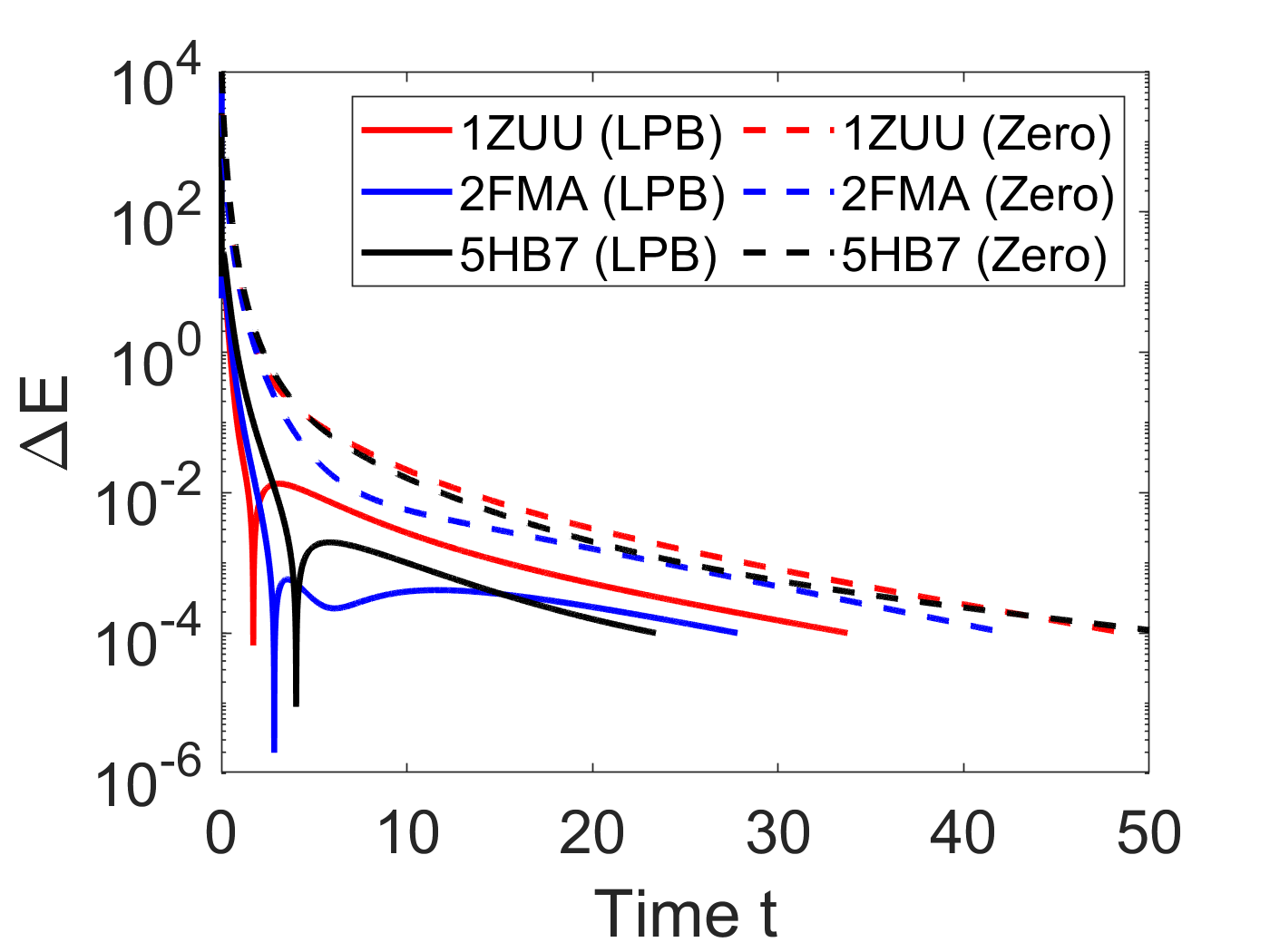}
    \caption{Constant $\dt$}
    \label{subfig:initialcond-de-const}
    \end{subfigure}
    \begin{subfigure}[b]{.45\textwidth}
    \includegraphics[width=\textwidth]{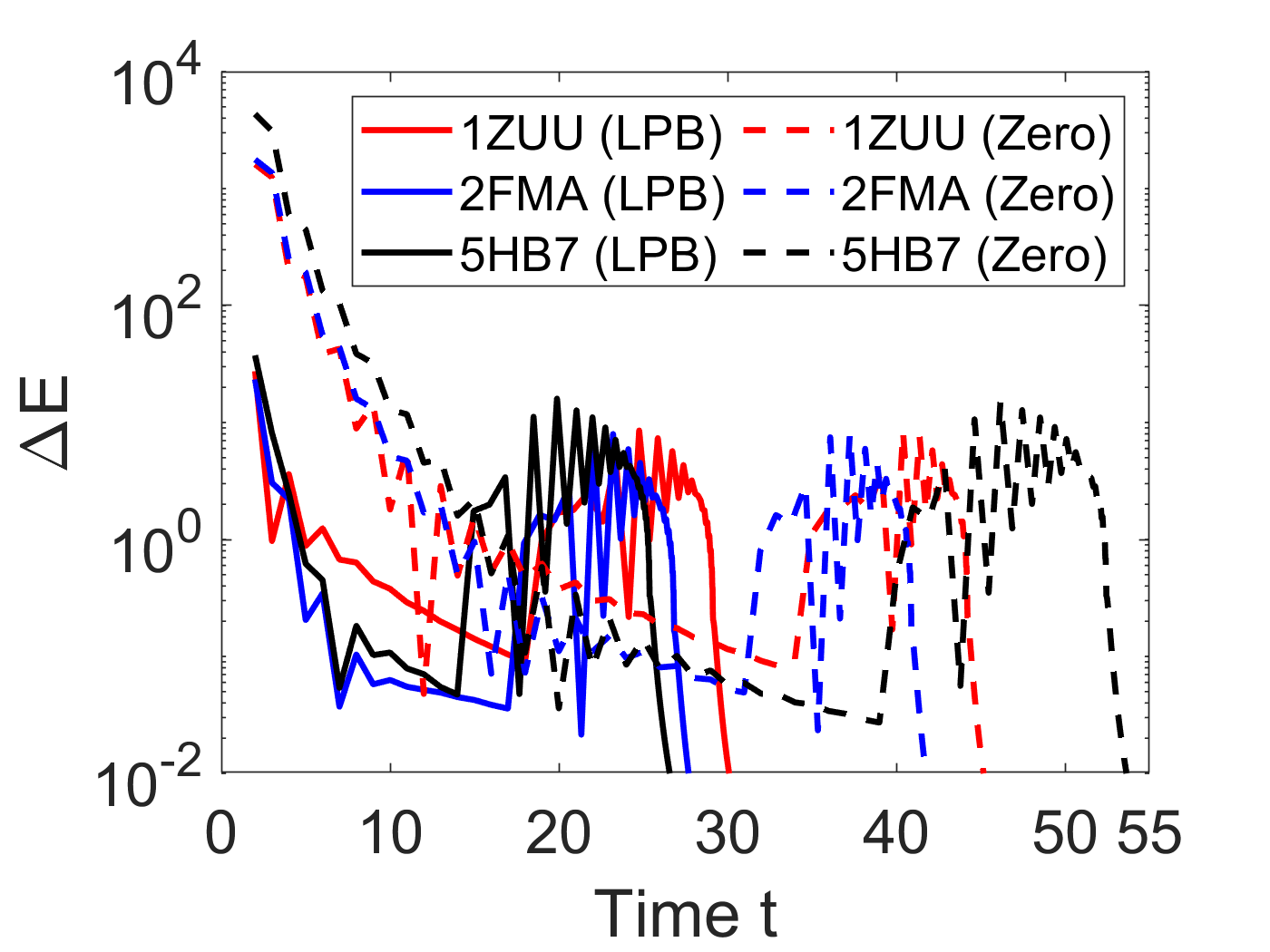}
    \caption{NonincreasingPID}
    \label{subfig:initialcond-de-nonincreasing}
    \end{subfigure}
    \caption{$\Delta E^n_{sol}$ vs time for proteins 1ZUU, 2FMA, and 5HB7 solved with the zero  and LPB initial conditions. The method for calculating  $\dt$ is constant in (a) and NonincreasingPID in (b).}
    \label{fig:initialcond-de}
\end{figure}

In Fig. \ref{fig:initialcond-de}, the energy difference $\Delta E^n_{sol}$ is plotted against the time $t$ for the constant $\dt$ and NonincreasingPID schemes in part (a) and (b), respectively. When using the LPB initial solution and holding $\dt$ as a constant, $\Delta E^n_{sol}$ for some proteins has an early, rapid dip that can be seen in Fig. \ref{fig:initialcond-de} (a). The three proteins, 1ZUU, 2FMA, and 5HB7, are selected because they all exhibited this early dip behavior. Since we use a tolerance condition on the energy difference, i.e., $\Delta E^n_{sol} < 10^{-4}$, as a stop condition, this would prematurely end the computation, resulting in a large error. To prevent this, we have added a condition to stop computation until after $t=5$, which passes this dip for the majority of proteins.  However, this does not catch all of them before the dip. Moreover, the choice of $t=5$ is entirely arbitrary, and there could be cases where this dip does not occur until much later. Therefore, in the present study when we generated the reference solution by using $\dt=0.01$ for 74 proteins, we always watch out the stopping time and
cross-validate the energies with those of the NoncreasingPID scheme. In case of a large disagreement, we re-generated the reference solution with $t=10$ or larger for preventing the early stop. 

As shown in Fig. \ref{fig:initialcond-de} (a), the free energies calculated based on zero initial condition is free of the early dip issue, and $\Delta E^n_{sol}$ decays monotonically as $t$ increases. So reference solution could also be generated by the constant $\dt$ method with zero initial condition. This is not adopted in the present study. First, the energies calculated in this way are essentially the same, as illustrated in Table \ref{tab:initialcond}. Second, as can be seen in both Fig. \ref{fig:initialcond-energy} and Fig. \ref{fig:initialcond-de}, with the zero initial condition, the constant $\dt$ computation will last longer. Consequently, it is more expensive to generate the reference solution by the zero initial condition than that by the LPB initial condition. 

The energy difference of the NonincreasingPID scheme shows some complicated patterns in Fig. \ref{fig:initialcond-de} (b). Recall that the stopping condition for the NonincreasingPID scheme is $\Delta E^n_{sol} < 10^{-2}$. This figure shows that the condition $\Delta E^n_{sol} < 10^{-1}$ may suffer an early stop problem, while $TOL=10^{-2}$ is completely safe. For both zero and LPB initial conditions, as $t$ starts from zero, $\Delta E^n_{sol}$ becomes smaller. Later, $\Delta E^n_{sol}$ rises up and involves some oscillations. During this period, $\dt$ is actually changed from $\dt_{max}$ to $\dt_{min}$. Once $\dt=\dt_{min}$, $\Delta E^n_{sol}$ shows an exponential decay and reaches the stopping condition $\Delta E^n_{sol} < 10^{-2}$ quickly. This study demonstrates that the NonincreasingPID scheme is a robust adaptive $\dt$ method. Moreover, the use of the zero initial condition in the NonincreasingPID scheme also needs a longer CPU time, such that $T_{end}=50$ is not enough for the protein 5HB7. In general, the NonincreasingPID scheme together with the LPB initial condition is what we recommended for future studies. 

\section{Conclusion}
In this paper, the pseudo-time GFM approach has been improved in two aspects. First, an Eulerian Solvent Excluded Surface (ESES) is implemented in our package to replace the MSMS for the molecular surface definition, which completely avoids the troublesome procedure of the Lagrangian-to-Eulerian conversion. The electrostatic analysis shows that the ESES free energy is more accurate than that of the MSMS, while being free of instabilities issues. The MSMS energy will converge to that of the ESES, as the MSMS density increases. However, for real proteins, the use of a high density usually induces instability in the PBE simulation. Second, various adaptive time integrations have been examined for the pseudo-time GFM approach. The usual hypothesis for steady state convergence assumes that one could use a longer $\dt$ as $t$ increases to save the CPU time. However, this hypothesis is found to be invalid for the pseudo-time GFM approach, perhaps because of the PBE nonlinearity and its time splitting treatment. A robust NonincreasingPID scheme is designed by tracking the temporal changes by a PID procedure, while maintaining the monotonic decay of $\dt$ with respect to $t$. In a large scale study with 74 proteins, the NonincreasingPID scheme yields an extremely accurate and efficient solution with a mean relative error of $0.06\%$ and an average saving  of $83.82\%$ CPU time. 

The adaptive time integration techniques developed in this work could be applied to other ADI algorithms \cite{Zhao15,Li20} for solving parabolic interface problems. In the future, we plan to explore more advanced adaptive time stepping methods for the pseudo-time PBE simulations.  A space dependent $\dt$ will be experimented to see if the steady state could be reached sooner.

\bigskip
\noindent
{\bf Acknowledgments}\\
This research is partially supported by the National Science Foundation (NSF) of USA under grant DMS-1812930.

\bibliographystyle{abbrv}
\bibliography{ESESDT-3}

\end{document}